\input jmacros.tex
\input epsf

\phantom{a}
\vskip.3in

{{\small \centerline{{\bf{THE HOMFLY POLYNOMIAL FOR LINKS}}}}
{{\small \centerline{{\bf IN RATIONAL HOMOLOGY 3-SPHERES}}}}
\bigskip
{\baselineskip=10pt{\msmall
\centerline{{\small E}{\msmall FSTRATIA } {\small K}{\msmall ALFAGIANNI\ \
{\msmall AND}\ \  {\small X}{\msmall IAO}-{\small S}{\msmall ONG}
{\small L}{\msmall IN }}}}
\bigskip
\centerline {\it September, 1995}
\bigskip
\p {{\small A}{\Small BSTRACT.}} {{ \msmall {
We construct a polynomial invariant, for
links in a large class of rational homology 3-spheres, which generalizes the
2-variable Jones polynomial (HOMFLY). As a consequence, we show that the
dual of the HOMFLY skein module of a homotopy 3-sphere is isomorpic to that 
of the genuine 3-sphere
.}}}}
\vskip.5in

{\baselineskip=11pt
\footnote { }{\p{\msmall The second author is partially supported by NSF.}}}
 
\p{\bf \S 1. Introduction}}
\medskip
The theory of quantum groups
gives a systematic way of producing families
of polynomial invariants, for knots and links
in $\R^3$ or $S^3$ (see for example [17]). In particular, the Jones polynomial
[9] and its generalizations ([5,10]),
can be obtained that way. All these Jones-type
invariants are defined as state models on knot
diagrams or as traces of braid group representations,
and the proofs of their topological invariance
offer little insight into the underlying topology.
The lack of a topological interpretation of these polynomial
invariants always makes it hard, if not impossible, to generalize them
for knots
and links in other 3-manifolds.

In his study of the topology
of the space of knots in $S^3$ [18], Vassiliev
constructed a vast family of knot invariants
which are now known as Vassiliev invariants or invariants
of finite type. Some of the nice aspects of the theory
of invariants of finite type are that they have a simple combinatorial 
description
and that they provide a unifying way to view
various knot polynomials ([2,3,13]). On the other hand, the
fact that the theory  of finite type invariants rests on topological
foundations allows its generalization to knots
in arbitrary 3-manifolds ([12, 15]). In this paper, using the machinery
developed in [12,15], we will show the existence and
uniqueness of link polynomials obeying the HOMFLY
skein relation in a  large class of rational homology 3-spheres.

In order to state  our main results we need to introduce some notation:
Suppose that $M$ is an orientable 3-manifold. Let $\pi=\pi_1(M)$
and let ${\hat \pi}$ denote the set of {\it non-trivial} conjugancy classes of $\pi$.
Notice that $\hat \pi$ can be identified with the set of 
non-trivial free homotopy classes of oriented loops in $M$. 
An $n$-component  {\it link} is a collection of $n$ unordered
oriented circles, tamely and disjointly embedded in $M$.
Hence, a link is homotopically equivalent to an unordered $n$-tuple
of elements in $\hat\pi\cup\{1\}$. In every homotopy
class of links, we will fix, once and for all, a link $CL$ and call it a
{\it trivial link}.
If $CL$ has $k$ components which are homotopically
trivial, our choice will be such that $CL=CL^*\coprod U^k$, where $U^k$
is the standard unlink with $k$ components in a small ball neighborhood
disjoint from $CL^*$ and $U^1$ will be abbreviated to $U$ later on. 
We will denote by $\cal C \cal L^*$ the set of all
{\it trivial
links} with none of their components homotopically trivial. It is obvious that
$\cal C \cal L^*$ is in 1-1 correspondence with unordered $n$-tuples of 
elements
in $\hat\pi$, for all $n>0$.

Every link $L$ is homotopic to a certain $CL^*\coprod U^k$ for some $CL^*\in
\cal C \cal L^*$, possibly empty. But the aim of link theory in the 3-manifold
$M$ is to understand how two links can differ up to a (tame) isotopy 
if they are
homotopic. Let $\cal L$ be the set of isotopy classes of links in $M$
and let $R=\C[v^{\pm1},z^{\pm1}]$ be the ring of Laurent polynomials in $v$
and $z$. A map ${\cal L}\rightarrow R$
will be called
a {\it link polynomial}.

Now let $z=t^{1\over 2}-t^{-{1\over 2}}$ and let $\cal I$ be the ideal of $R[t]$ generated by
$v-v^{-1}$ and $t$. Let $\hat R$ be the pro-$\cal I$ completion of $R[t]$, i.e. the
inverse limit of 
$$\cdots\rightarrow R[t]/{\cal I}^n\rightarrow R[t]/{\cal I}^{n-1}\rightarrow\cdots.$$

\medskip
\p {\bf Theorem A.}  {\sl Let $M$ be a rational homology 3-sphere
which is either atoroidal or a Seifert fibered space. Then, there is a unique
map $J_M: {\cal L}\rightarrow\hat R$ satisfying the HOMFLY
skein relation
$$v^{-1}J_M(L_+)-vJ_M(L_-)=zJ_M(L_o)$$
and with given values  $J_M(U)$ and $J_M(CL^*)$ for every $CL^*\in\cal C\cal L^*$.
Moreover, we may choose $J_M(U)$ and $J_M(CL^*)$ in $R$ appropriately such that
$J_M$ is a link polynomial, i.e. $J_M(L)\in R$ for every $L\in\cal L$.}
\medskip

Here, as usual, the three links $L_+$, $L_-$ and $L_o$ appeared in the HOMFLY
skein relation differ only in a small ball neighborhood in $M$ where, under a
suitable projection, they intersect at a positive crossing, a negative 
crossing, and a smoothing of a crossing, respectively.

It seems to be worthwhile to comment on the relationship between Theorem A and
the study of the HOMFLY skein modules of 3-manifolds (see [7] for a survey).
With the notation as above, the HOMFLY skein module ${\cal S}_3(M)$ is defined
to be the $R$-module spanned by $\cal L$, and subject to the HOMFLY skein
relation
$$v^{-1}L_+-vL_-=zL_o.$$

Let $S(R\hat\pi)$ be the symmetric tensor algebra
of the free $R$-module
$R\hat\pi$ generated by $\hat\pi$. J. Przytyski proposed
the following conjecture.

\medskip
\p {\bf Conjecture.} {\sl If $M$ is compact and contains no 
closed non-separating
surfaces, then
$${\cal S}_3(M)\cong S(R\hat\pi)$$
as $R$-modules.}
\medskip

In the special case of $M=X\times [0,1]$ where $X$ is a compact surface, it was
proved by Przytyski (see [7] and references therein) that  ${\cal S}_3(M)\cong
S(R\hat\pi)$ as $R$-algebras. Many other partial verifications of the conjecture
are known. In the
presence of non-separating closed surfaces in $M$, one can construct
examples in which ${\cal S}_3(M)$ is not torsion free.

Notice that ${\cal CL}^*\cup\{U\}$ is in one-one correspondance with a basis of
$S(R\hat\pi)$. So Theorem A can be thought of providing a strong supporting evidence
to (the dual version) of Przytyski's conjecture.  In the special case when $\pi_1(M)=1$,
Theorems A is of particular
interest because of the Poincare conjecture. Also, in this special case, we may  rephrase
Theorem A in the language of skein modules since ${\cal CL}^*=\emptyset$ when
$\pi_1(M)=1$. 

We denote by ${\cal S}_3^*(M)$ the dual of ${\cal S}_3(M)$, i.e. ${\cal S}_3^*=\hbox{Hom}
({\cal S}_3(M),R)$.

\medskip
\p {\bf Theorem B.} {\sl If $M$ is a homotopy 3-sphere, then
$${\cal S}_3^*(M)\cong{\cal S}_3^*(S^3)\cong R$$
as $R$-modules.}
\medskip

Let us now briefly discuss the main ideas of the paper:

For an $L\in {\cal L}$ let $ {\cal M}^L$ denote the space
of maps, equipped with the compact-open topology,
from a disjoint union of circles to $M$, which are homotopic to $L$.
Any two links in some $ {\cal M}^L$ are related
by a sequence of ``crossing changes". When we make a
crossing change from one link to another,
we produce a singular link with one transverse double point
as an intermediate step. If we are given a
link invariant $F: {\cal L}\longrightarrow
{\cal R}$, where $\cal R$ is a ring, we may define an invariant $f$ 
of singular 
links with one double point by 
$$f(L_{\times})=F(L_+)-F(L_-),\eqno (2)$$
where $\times$ stands for a double point and $L_\pm$ are links obtained from
resolving the double point into a positive or negative crossing, 
respectively.

A natural question is the following: Starting with a singular link 
invariant $f$,
we want to find necessary and sufficient conditions
that $f$ has to satisfy so that it is derived from
a link invariant, via (2). This question and its general version
was shown in [12] and [15] to play an important role in understanding 
invariants
of finite type for links in 3-manifolds.
In $\S 3$ we answer the question for links in
rational homology spheres which are either {\it atoroidal}
or {\it Seifert fibered spaces}. In \S 4 and \S 5, we use the results
in $\S 3$
to prove Theorem A. Theorem B is a direct corollary of Theorem A.

The question of whether $f$ is induced by a link invariant
turns out to be a question about
the {\it integrability} of $f$ along
paths in every ${\cal M}^L$. Let
$\Phi$ be a path in ${\cal M}^L$ connecting two links. After
perturbation,
we may assume that there are only finitely many points on $\Phi$ 
where we see singular links with one double point.
Moreover, we can assume that
when the parameter of the path passes through a point
where we see a singular link, the nearby links are changed by
a crossing change. 
The sum of
suitably signed values of $f$ on these singular links along the 
path $\Phi$, 
denoted here by $X_{\Phi}$,
can be
thought of as the integral of $f$ along $\Phi$.

In order for $f$ to be derived from a link
invariant, it is necessary and sufficient
that  $X_{\Phi}$ is independent of $\Phi$ relative to the end points, or
equivalently that $X_{\Phi}=0$, for every loop
$\Phi$ in ${\cal M}^L$. It turns
out that the first thing that one has to do
is to find a set
of finitely many {\it local integrability conditions}
which guarantee $X_{\Phi}=0$ for every null-homotopic loop $\Phi$
in $\cal M^L$. Our technical assumptions about 3-manifolds will then 
imply that these local integrability conditions guarantee $X_{\Phi}
=0$ for every loop $\Phi$ in $\cal M^L$.

If $\Phi$ is null-homotopic in ${\cal M}^L$ then we achieve our
goal by putting the null-homotopy into {\it almost general position}. This
is done in [15] (see $\S 3.2$). The treatment for the general
case is based on the machinery developed in [12]. Here we need to change our
point of view and think of $\Phi$ as a map from a disjoint union
of tori into  $M$. This naturally leads to the study
of tori in  $M$ and to the use of the results in [8] and [16] in
order to treat essential $\Phi$'s.
More precisely, by employing
the homotopy classification of essential tori in 
Seifert fibered spaces, we are able to homotope $\Phi$ into
a certain nice position so that the local
integrable conditions imply $X_{\Phi}=0$ (see $\S 3.3$ and $\S 3.4$). 

We organize the paper as follows:
In $\S 2$ we recall the generic picture of a family of maps
from a compact 1-polyhedron into a 3-manifold, parameterized by a disc,
and we give the preliminaries from the topology of 3-manifolds that we
use in subsequent sections. In $\S 3$ we answer the integrability question 
addressed above.
The main result of this section is Theorem 3.1.2.
In $\S 4$ we use Theorem 3.1.2 to construct
formal power series link invariants which satisfy the same crossing 
change formulae as the HOMFLY polynomial (Theorem 4.2.1).
In $\S 5$ we  use the surgery description of 3-manifolds
to prove that, with an appropriate choice of the initial values, the power series obtained
in $\S 4$ converge to Laurent polynomials. Our main result is Theorem 5.2.3. 
By combining
Theorem 4.2.1 and Theorem 5.2.3, Theorem A, and subsequently Theorems B, will follow
immediately.
\medskip

\noindent{\sl Acknowledgment.} We would like to thank Joan Birman
for many helpful conversations during the early stages
of this work. We also thank
Peter Scott for helpful correspondence on the
topology of 3-manifolds and for his help with
the proof of Lemma 3.3.2. Finally, this work was initiated and completed during
the second and the first authors' visits (1993-94 and 1994-95, respectively) to the
Institute for Advanced Study. We would
like to thank IAS for its hospitality. 
\bigskip
%\vfill
%\eject

\p{\bf \S 2. Preliminaries}
\medskip
\p {\bf \S 2.1. Almost general position for a disjoint union of circles}
\medskip

In this section we summarize from [15] the results about
the generic picture of a family of maps from
a disjoint union of circles to a 3-manifold,
parameterized by a disc.

Let $P$ be an 1-dimensional compact polyhedron. Let $M$ be a 3-manifold and
let
$D^2$ denote the 2-disc.
A map $\Phi :P\times D^2\longrightarrow M$
gives rise to a family of maps $\{\phi _x : P\longrightarrow M
\ ; \ x \in D^2\}$ where $\phi_x(*)=\Phi(*, x)$ for $x\in D^2$.
Suppose that every $\phi_x$ is a piecewise--linear map and let
$S_{\phi}$ be the closure of the set
$\{ x\in D^2 \ ; \ \phi_x \hbox { is not an embedding }\}$.
One can see that $S_{\phi}$ is a sub--polyhedron of $F$.

Two maps $\phi_1, \phi_2 : P\longrightarrow M$ are called
ambient isotopic if there exists an isotopy $h_t:
M\longrightarrow M$, $t\in [0, 1]$ with $h_0=id$ and
$h_1\phi_1=\phi_2$.

Let us now introduce some terminology about $1$--dimensional polyhedra in
3- manifolds.

Let $P$ be an $1$--dimensional polyhedron. Every point
$q\in P$ has a neighborhood homeomorphic to a bouquet
of finitely many arcs such that $q$ is the common endpoint of
these arcs. The number of arcs in the bouquet is called the
{\it valence} of $q$. A point $q\in P$ with valence different than
$2$ is called a {\sl vertex} of $P$.
A component of the complement of vertices is called an
{\sl edge} of $P$.

A double point  of a map
$\phi:P\longrightarrow M$ is a point $p\in M$ such that
$\phi^{-1}(p)$ consists of two points. A
double  point of a piecewise linear map $\phi :
P\longrightarrow M$ is called transverse double point
 if there exist two
$1$--simplexes $\sigma_1$, $\sigma_2$ contained in the $1$--skeleton of $P$ such
 that

 1)  $\phi$ is linear and non--degenerate on $\sigma_1$ and
$\sigma_2$ ;

 2) $\phi(\sigma_1)\cap \phi(\sigma_2)$
is the double in question;

 3) $\phi(\sigma_1)$ and $\phi(\sigma_2)$ intersect
transversally in their interiors and they
do not lie on the same plane.

We call an $1$--dimensional subpolyhedron
$S\subset D^2$ {\sl neat} if $S\cap \partial D^2$ consists of finitely
many points and each of the them is a valence $1$ vertex of
$S$. We call these vertices {\sl boundary vertices} of $S$ and
we call the vertices of $S$ lying in the interior of $D^2$
{\sl interior vertices} of $S$.
\smallskip
We suppose now that $P$ is a disjoint union of circles. Then we have

\medskip
\p {\bf Proposition 2.1.1.} ([15]) {\sl A map $\Phi :
P\times D^2 \longrightarrow M$ can be changed by an arbitrary
small perturbation so that $S_{\phi}$ is a neat
$1$--dimensional subpolyhedron of $F$. Moreover, we have

1) If $x, x'\in D^2$ belong to the same component of
$D^2\setminus S_\phi$ or $S_\phi\setminus$ ${\{}$ interior
vertices$\}$, then $\phi_x$ and $\phi_{x'}$ are ambient
isotopic;

2) the interior vertices of $S_\phi$ are of valence either
four or one;

3) If $x\in S_\phi$ lies on an edge of $S_\phi$ or is a
boundary vertex, then $\phi_x$ has exactly one transverse
double point;

4) If $x\in S_\phi$ is an interior vertex of valence four,
then $\phi_x$ has exactly two transverse double points;
and

5) If $x\in S_\phi$ is an interior vertex of valence one,
then $\phi_x$ is an embedding ambient isotopic to the nearby
embeddings.}
\medskip

We say that the resulting map in Proposition 2.1.1 is
{\sl in almost general position}. Figure 2.1 below
illustrates $S_{\Phi}\subset D^2$ for a
map $\Phi$ in almost general position.

\bigskip
\centerline{\epsfysize=1.25in\epsfbox[49 56
566 735]{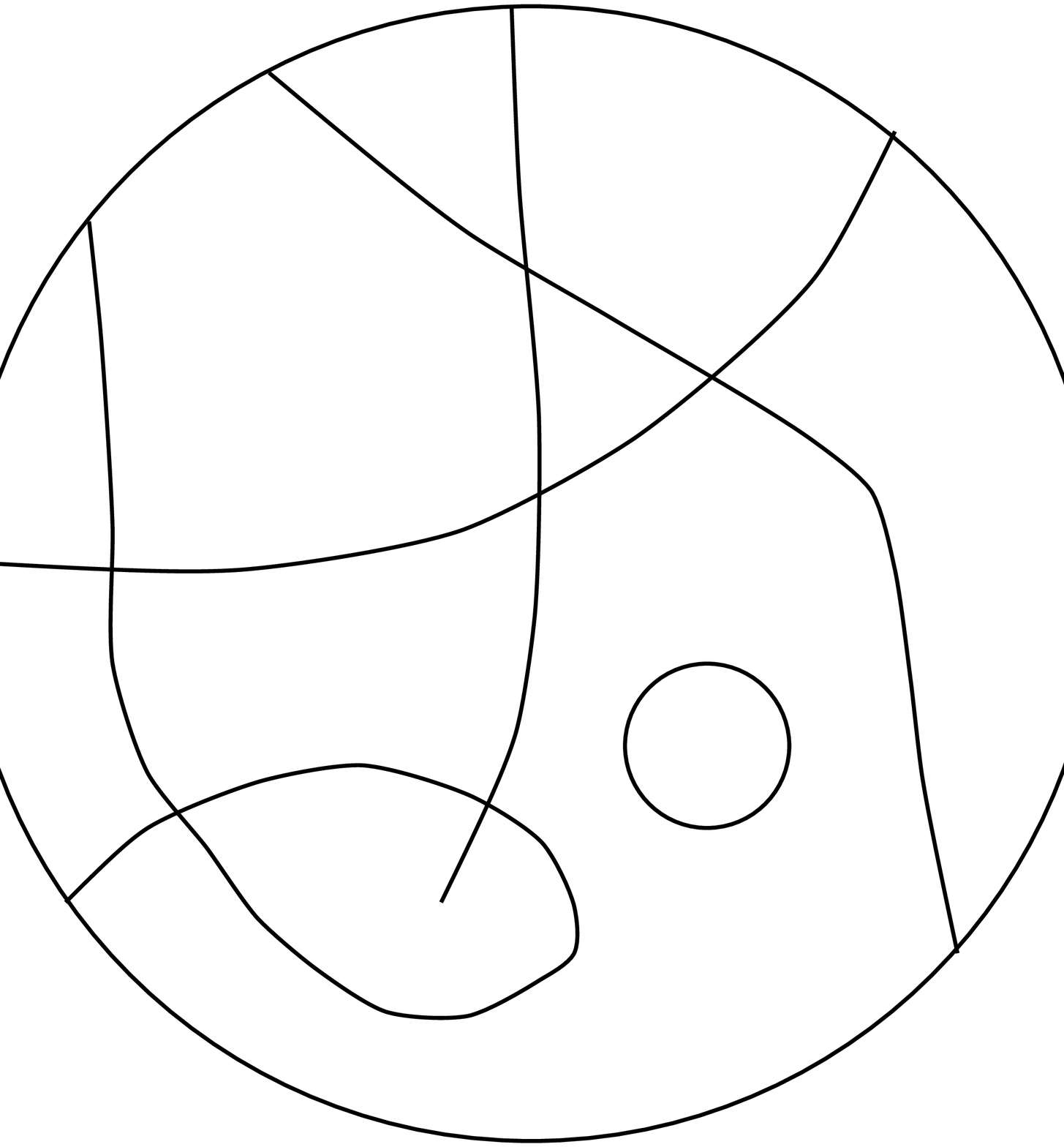} }
\smallskip
\centerline {$S_{\Phi} \subset D^2$.}
\centerline {{\bf Figure 2.1}}
\medskip

\medskip
\p {\bf Remark 2.1.2.}  a) The proposition above is true for
maps $\Phi : P\times X \longrightarrow M$, where $X$ is a compact surface with 
$\partial
 X \neq \emptyset$, as well as more general types
of polyhedra $P$. It was used extensively, in [12] to define and study
invariants of {\it finite type} for knots in 3-manifolds.

\p b) If ${\Phi} | {P\times {\partial D^2}}$, is in almost general position 
already, then the perturbation in proposition 2.1.1 can be made
relatively $\partial D^2$.
\medskip

\medskip
\p {\bf \S 2.2. Seifert fibered spaces}
\medskip

In this section we give some terminology from the topology of 3-manifolds 
and recall some results from  [8] and [16] that are used in 
subsequent sections.

\medskip
\p {\bf Definition 2.2.1.} {\sl A surface $X\ne S^2$, properly embedded
in a 3-manifold $M$ (or embedded in $\partial M$), is compressible
if there exists a disc $D\subset M$ such that $D\cap X= \partial D$
and $\partial D$ is not homotopically trivial in $X$. Otherwise
$X$ is called incompressible in $M$.
A compact, orientable, irreducible 3-manifold is called a {\it Haken}
(or {\it sufficiently large}) manifold, if it contains a two-sided
incompressible surface.}
\medskip

\medskip
\p {\bf Definition 2.2.2.} {\sl Let $M$ be a closed 
3-manifold and $X$ a surface. A map $\phi: X\longrightarrow M$ is called
essential if $\hbox{ker}\{\phi_*: \pi_1(X) \longrightarrow
\pi_1(M)\}=1$.}
\medskip

Let $(\mu,\nu)$ be a pair of relatively prime integers. Let
$$ D^2= \{ (r, \ \theta); \ \ 0\leq r \leq 1,\ 0\leq \theta
\leq 2\pi \}\subset {\R}^2$$
\p  A fibered solid torus
of type $(\mu, \ \nu)$, is the quotient of the cylinder $D^2\times I$
via the identification $((r, \ \theta),\ 1)=((r, \ \theta+
{{2\pi \nu }\over {\mu}}), \ 0)$. If $\mu > 1$ the fibered solid
torus is said to be {\it exceptionally fibered} and the core is the
{\it exceptional fiber}. Otherwise the fibered solid torus is {\it regularly
fibered} and each fiber is a {\it regular fiber}.

\medskip
\p {\bf Definition 2.2.3.} {\sl An orientable 3-manifold is said to be a
{\it Seifert fibered space}, if it is a union of pairwise disjoint simple
closed curves, called {\it fibers}, such that each one has a closed
neighborhood, consisting of a union of fibers, which is homeomorphic
to a fibered solid torus via a fiber preserving isomorphism.}
\medskip

In a Seifert fibered space $M$, a fiber is called {\it exceptional} if it has a
neighborhood homeomorphic to an exceptionally fibered solid torus and the fiber
in question corresponds to the exceptional fiber
of the solid torus.
The {\it orbit space} of $M$ is the quotient obtained by identifying
every fiber of $M$ to a point.

\medskip
\p {\bf Definition 2.2.4.} {\sl Let $M$ be a Seifert fibered space, with
a fixed fibration and let  $\break {p: M\longrightarrow B}$
be the fiber projection.
Let $X$ be a surface. A map $\phi: X\longrightarrow M$ is called
{\it vertical} or {\it saturated}, with respect to $p$,
if $p^{-1}(p\phi(X))=\phi(X)$ and $\phi(X)$ contains no
exceptional fibers.}
\medskip

Throughout this paper we are dealing with Seifert fibered spaces
that are rational homology 3-spheres (i.e $H_1(M)$ is finite).
It is known that if $\pi_1(M)$ is infinite then either $M$ is Haken
or it fibers over the 2-sphere with three 
exceptional fibers. 
The following  proposition gives a classification up to homotopy of singular
essential tori in such manifolds.

\medskip
\p {\bf Proposition 2.2.5.} { \sl Suppose that $M$ is a Seifert fibered
rational
homology 3-sphere with infinite $\pi_1$. Let
$\Phi: T=\coprod T_i \longrightarrow M$ be an essential map,
where each $T_i$ is a torus.
Then there exists a homotopy
$\Phi_t: T\longrightarrow S$, $t\in [0, \ 1]$, such that
$\Phi_0=\Phi$ and $\Phi_1$ is vertical with respect
to the fibration of $M$.}
\medskip

The proof of the proposition is given in Proposition 5.13 of [8]
for the case when $M$ is Haken. For the case that $M$ is not Haken,
see Theorem 6.4 of [16].
\bigskip

\p {\bf \S 3. Integrating invariants of singular links}
\medskip

In this section we introduce singular links and study their invariants.
Our purpose is to give conditions under which an invariant of singular links
gives rise to a link invariant.

\medskip
\p {\bf \S3.1. Definitions and the statement of the main result}
\medskip

Let $M$ be an oriented 3-manifold
and  let $P$ be a disjoint union of oriented circles.

\medskip
\p {\bf Definition 3.1.1.} {\sl A singular link of order $n$ is a
piecewise-linear map $ L: P \longrightarrow M$ that
has exactly $n$ transverse double points. Two singular links $L$ and $L'$
are equivalent if there is an isotopy $ h_t: M \longrightarrow M$, $
t\in [0,1]$ such that $h_0=id$, $L'=h_1(L)$ and the double points of
$h_t(L)$ are transverse for every $t\in [0,1]$}.
\medskip

We will also use $L$ to denote $L(P)$. A singular link of order $0$ is
simply a link.

Let $p\in M$ be a transverse double point of a singular link $L$.
Then $L^{-1}(p)$ consists of two points $p_1, \ p_2 \in L$. There are
disjoint 1-simplexes, $\sigma_1$ and $\sigma_2$, on $P$  with $p_i
\in\hbox{int}(
\sigma_i)$, $i=1, \ 2$ such that for a small ball neighborhood $B$ of $p$
in $M$
$$ L \cap B= L(\sigma_1)\cap L(\sigma_2)$$
Moreover, there is a proper 2-disc $D$ in $B$ such that $L(\sigma_1)$, $L(\sigma
_2) \subset D$ intersect transversally at p, and the isotopy $h_t$ of
definition 3.1.1 carries the ball disc pair $(B ,  D)$ through for all the
double points of $L$.

We can resolve a transverse double point of a singular knot of order $n$
in different ways. Notice that $L(\sigma_1)\cap L(\sigma_2)$ consists
of four points on $\partial D$. Also, since $\sigma_i$
inherits an orientation from that of $P$ we can talk about the initial point
and terminal point of $\sigma_i$ and $L(\sigma_i)$.

Now choose arcs $a_1$, $a_2$,
$b_1$, $b_2$ with disjoint interiors such that

 (1) $a_1$ and $a_2$ go from the initial point of $L(\sigma_1)$
to the terminal point of $L(\sigma_1)$ and lie in distinct components
of $\partial B \setminus \partial D$; and
 
 (2) $b_1$ and $b_2$ lie on $\partial D$ with $\beta_1$ going from
the initial point of $L(\sigma_1)$
to the terminal point of $L(\sigma_2)$ and $b_2$
from the initial point of $L(\sigma_2)$ to the terminal point of
$L(\sigma_1)$. See Figure 3.1.

\bigskip
\centerline{\epsfysize=1.5in\epsfbox[49 56 566 735]{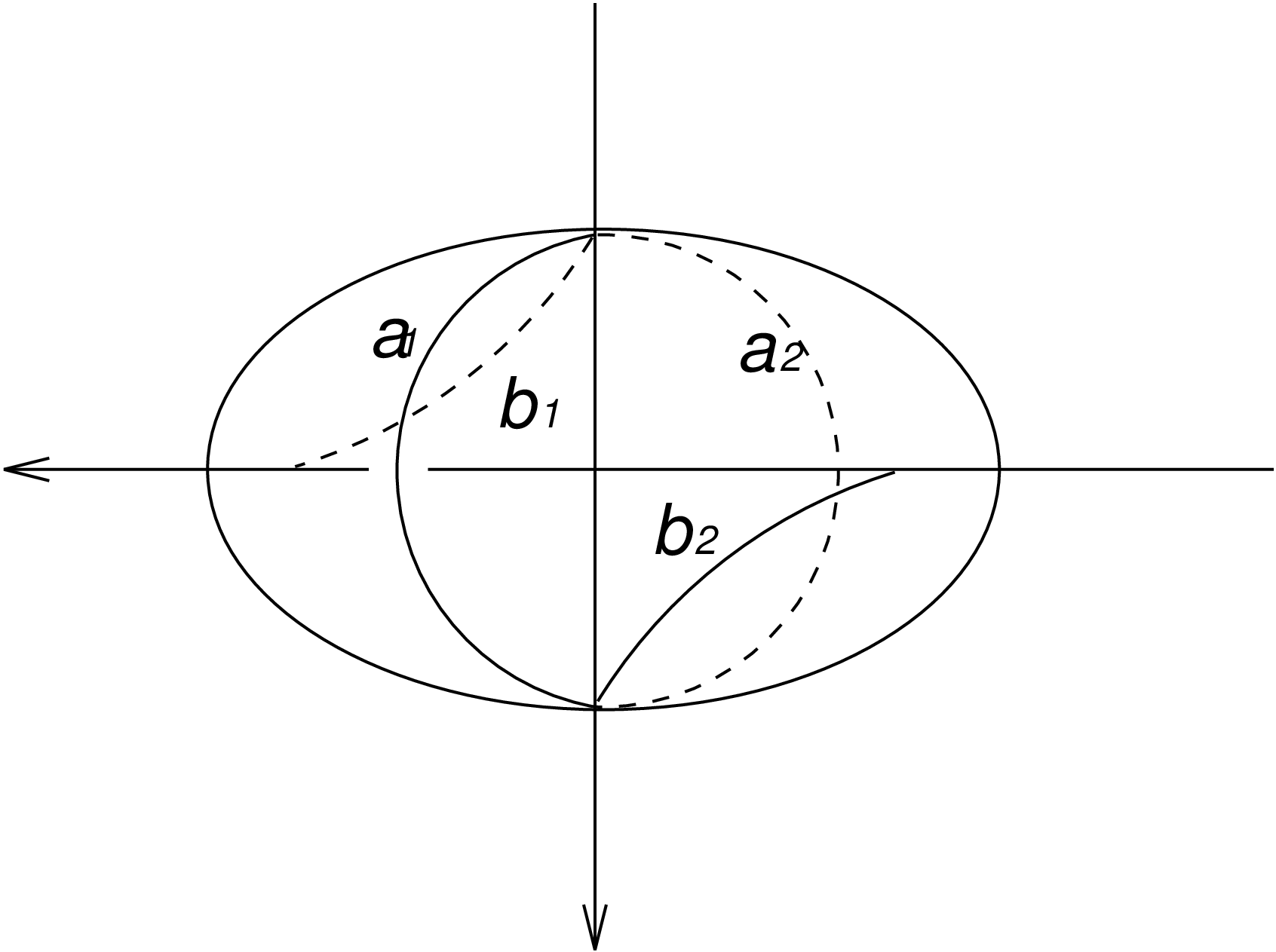} }
\smallskip
\centerline {Resolving a transverse double point.}
\centerline {{\bf Figure 3.1}}
\medskip

Define
$$L_+=\overline {K\setminus K(\sigma_2)} \cup a_1$$
$$L_-=\overline {K\setminus K(\sigma_2)} \cup a_2$$
$$L_0=\overline {K\setminus K(\sigma_2) }\cup (b_1\cup b_2)$$

Clearly $L_+$, $L_-$ are well defined singular links
of order ${n-1}$. We call $L_+$ (respectively $L_-$) a positive
(respectively a negative) resolution of $L$.

We will denote by ${\cal L}^{n}$ the set of equivalence classes of
singular
links of order $n$ in $M$. Let $\cal R$ be a ring. A singular link invariant
is a map ${\cal L}^{n} \longrightarrow \cal R$. Notice that for $n=0$
we have a link invariant. From a link invariant
$F:{\cal L} \longrightarrow \cal R$ we can always define a singular
link invariant $f:{\cal L}^{1} \longrightarrow \cal R$ as follows:

Let $L_{\times} \in {\cal L}^{1}$ where $\times $ stands
for the only double point. Then $L_+$, $L_- \in {\cal L}^{0}=
{\cal L}$. We define $ f: \ {\cal L}^1\longrightarrow \cal R$ by
$$f(L_{\times})= F(L_+)-F(L_-) \eqno (3.1.1)$$

In this section we will answer the following question:
Suppose that we are given a singular link invariant
$f: {\cal L}^1\longrightarrow \cal R$. Under what conditions can
we find a link invariant $F:{\cal L} \longrightarrow \cal R$ so that
(3.1.1) holds for all $L_{\times} \in {\cal L}^{1}$.

In [2], D.
Bar-Natan thinks of $(3.1.1)$ when going from the link invariant $F$ to 
the singular link invariant $f$ as 
the ``first derivative"
of $F$. In this spirit
the question above concerns the ``integrability" of a singular link invariant.

For the rest of the paper we will assume, unless otherwise stated, that $M$ is
a rational homology 3-sphere such that either i) it has 
trivial $\pi_2$ and there 
are no essential maps $S^1\times S^1
 \longrightarrow M$ or
ii) it is a Seifert fibered space. If $M$ is as in i) we will say that
it is {\it atoroidal}. Notice that if $M$ is as in ii)
then it is {\it irreducible} and hence we have $\pi_2(M)=1$ 
by the sphere theorem (see for example [7]).

We will also assume that
$\cal R$ is a ring which is torsion free as an abelian group.
Our main result in this section is the following theorem,
which answers
the integrability question for a large class of rational homology 3-spheres.

\medskip
\p{\bf Theorem 3.1.2.} {\sl Suppose that $M$ and $\cal R$ are as above, and
let $f: {\cal L}^{1} \longrightarrow {\cal R}$
be a singular link invariant. There exists a link invariant
$F: {\cal L} \longrightarrow {\cal R}$ so that (3.1.1)
holds for all $L_{\times} \in {\cal L}^1$ if an only if
$f$ satisfies}
$$f({\epsfysize=0.075in\epsfbox{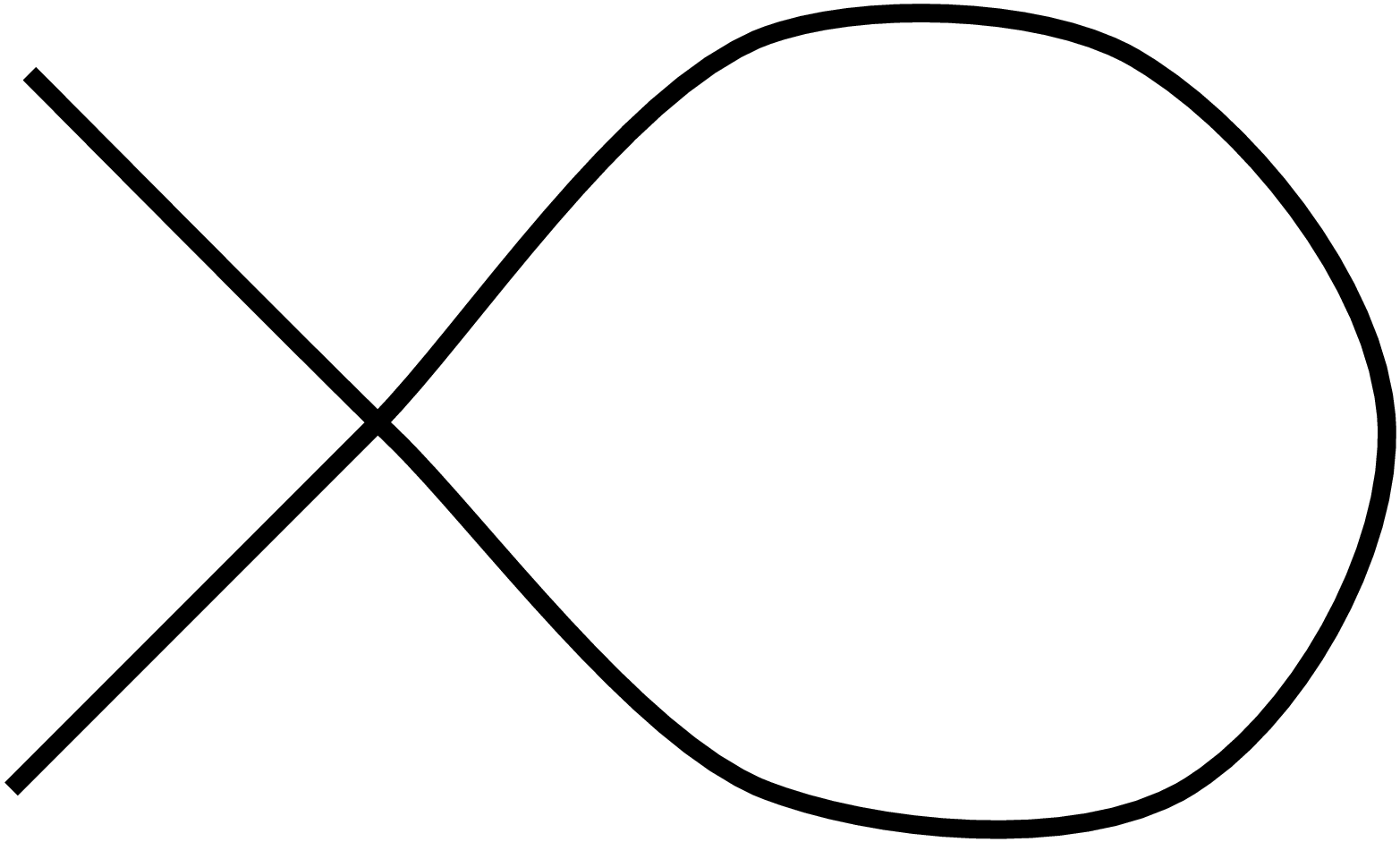} })=0  \eqno (1)$$
$$f(L_{{\times}+})-f(L_{{\times}-})=f(L_{+{\times}})-f(L_{-{\times}}) \eqno (2)
$$
\medskip

\medskip
\p{\it Notation:} Before we proceed with the proof of Theorem 3.1.2, let us 
explain the notation above. In (1) the kink stands for a singular link
$L_{\times} \in {\cal L}^{1}$ where there is a 2-disc $D\subset M$
such that $L_{\times} \cap D=\partial D$, and the unique double point
of $L_{\times}$ lies on $\partial D$. In (2) we start with an
arbitrary singular link $L_{\times \times} \in {\cal L}^{2}$. The
four singular links in ${\cal L}^{1}$ are obtained by
resolving one double point of $L_{\times \times}$ at a time.
Following [15], we will call conditions (1) and (2) above the {\it local
integrability conditions}.
\medskip

The proof of Theorem 3.1.2 will occupy the rest of \S 3.

\medskip
\p {\bf Proof of Theorem 3.1.2.} One direction of the theorem is clear. That is
if a singular
link invariant $f: {\cal L}^{1} \longrightarrow {\cal R}$ is
derived from a link invariant $F: {\cal L} \longrightarrow {\cal R}$
via (3.1.1), then it satisfies
(1) and (2). To see that (1) is satisfied observe that
the positive and the negative resolution of the double point in the kink
are equivalent. For (2) observe that, using (3.1.1),
both sides of (2) can be
expressed as $F(L_{++})-F(L_{-+})-F(L_{+-})+F(L_{--})$.

We now turn into the proof of the other direction. Namely, 
assuming that a singular link
invariant $f: {\cal L}^1\longrightarrow {\cal R}$ is given satisfying
(1) and (2), we show that it can be derived from a link invariant
$F: {\cal L} \longrightarrow {\cal R}$ via (3.1.1), provided that
$M$ is as in the statement of Theorem 3.1.2. The proof will be broken 
into several steps.

Let $L\in {\cal L}$ be a link in $M$.  We also use
$L$ to denote a representative $L:P \longrightarrow
 M$, of $L$. Let ${\cal M}^L(P, M)$ denote the space of maps
$P\longrightarrow M$ homotopic to $L$, equipped with the compact-open
topology. For every $L' \in {\cal M}^L(P, M)$,
we choose a homotopy ${\phi}_t: P {\times}[0, 1]
\longrightarrow M$ such that ${\phi}_0=L'$ and ${\phi}_1=L$.
After a small perturbation, we can assume that for only finitely many points
$0<t_1<t_2<\cdots<t_n<1$, ${\phi}_t$ is not an embedding. Moreover, we can
assume that ${\phi}_{t_{i}}$, for $i=1,2,\dots,n$ are singular
links of order $1$ (i.e. ${\phi}_{t_{i}} \in {\cal L} ^{(1)} $).
For different $t's$ in an interval of $[0, \ 1]\setminus \{t_1, \ t_2, \
\dots, t_n\}$, the corresponding links are equivalent. When $t$ passes
through $t_i$, ${\phi}_t$ changes from one resolution of ${\phi}_{t_i}$
to another.

We define
$${F(L')=F(L)+\sum_{i=1}^{n}{\epsilon}_i f({\phi}_{t_i})}\eqno (3.1.2)$$
Here ${\epsilon}_i={\pm}1$ is determined as follows: If
${\phi}_{t_i+\delta}$, for $\delta >0$ sufficiently small, is a positive
resolution of ${\phi}_{t_i}$ then ${\epsilon}_i=1$. Otherwise
${\epsilon}_i={-1}$.

To prove that $F$ is well defined we have to show that modulo ``the
integration constant" $F(L)$, the definition of $F(L')$ by $(3.1.2)$
is independent of the choice of the homotopy. For this we consider a closed
homotopy $\Phi:P\times S^1\longrightarrow M$. After a small perturbation,
we can assume that there are only finitely many points
$x_1,x_2,\dots,x_n \in S^1$, ordered cyclicly according to the
orientation of $m$, so that ${\phi}_{x_i}\in {\cal L}^{1}$ and
$\phi_x$ is equivalent to $\phi_y$ for all $x_i<x,y<x_{i+1}$. To prove that $F$
is well defined
we need to show that
$${X_{\Phi}:=
\sum_{i=1}^{n}{\epsilon}_i f({\phi}_{t_i})=0}\eqno (3.1.3)$$
where ${\epsilon}_i={\pm}1$ is determined by the same rule as above.
 We will call $(3.1.3)$
the {\it global integrability condition} around $\Phi$.

The proof of (3.1.3), which will be broken into many steps, occupies the rest of
\S3.

\medskip
\p {\bf \S3.2. The proof of the global integrability condition in some special
cases}
\medskip

Assume that $M$ is an oriented 3-manifold, with $\pi_2(M)=1$, and
that $f: {\cal L}^1 \longrightarrow {\cal R}$
is a singular link invariant. Let $L: P \longrightarrow M$ be a link,
and recall that  ${\cal M}^L(P, M)$ denotes the space of maps
$P\longrightarrow M$ homotopic to $L$, equipped with the compact-open
topology. A closed homotopy $\Phi: P\times S^1 \longrightarrow M$
from $L$ to itself, can be viewed as a loop in
${\cal M}^L(P, M)$.

\medskip
\p {\bf Lemma 3.2.1.} {\sl Let $M$, $P$, and $\Phi$ be as above. Moreover, 
suppose that $\Phi$ can be extended to a map $\hat\Phi: P\times D^2
\longrightarrow M$, where $D^2$ is a 2-disc with $\partial D^2 = \{*\} \times
S^1$. Then, $\Phi$ satisfies the global integrability condition, 
i.e. $X_{\Phi}=0$.}
\medskip

\p {\bf Proof.} We perturb ${{\hat \Phi}}$ to an almost general position map
as in Proposition 2.1.1.
Then each edge of the set of singularities  $S_{{\hat {\Phi}}}$,
corresponds to a singular link of order 1. So by using the invariant
$f$ we can assign an element of $\cal R$ to every edge of $S_{{\hat {\Phi}}}$.
we can reduce the global integral condition around ${\hat {\Phi}}$,
to local integrable conditions
around each interior vertex in $S_{{\hat {\Phi}}}$.

\bigskip
\medskip
\centerline{\epsfysize=1.5in\epsfbox[49 56 566 735]{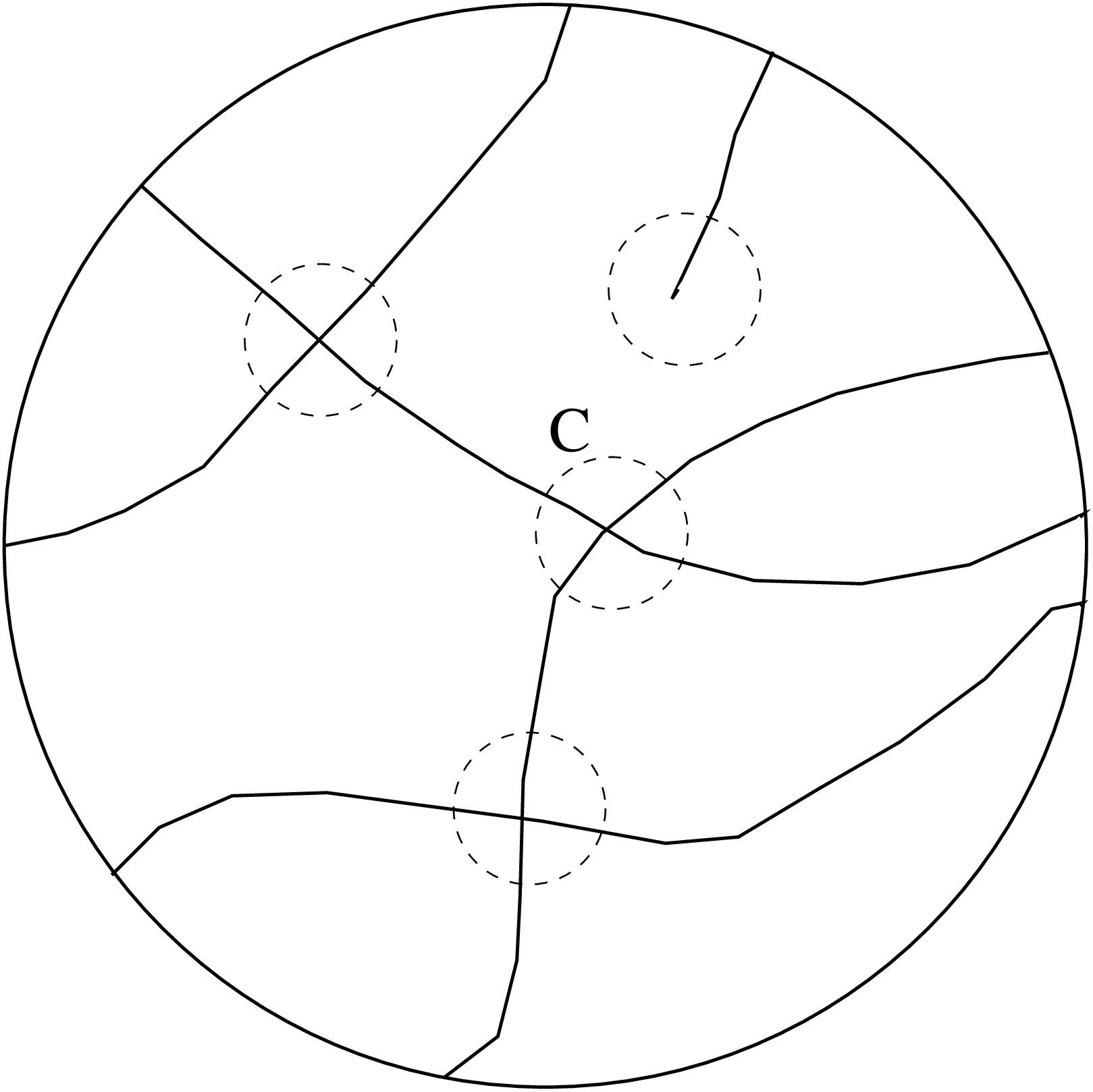} }
\bigskip
\centerline {From global to local integrability conditions.}
\centerline {{\bf Figure 3.2}}
\medskip

More precisely, for every interior
vertex of $S_{{\hat {\Phi}}}$ draw a small circle $C$ around it, so that the
number of
points in $C\cap S_{{\hat {\Phi}}}$ is equal to the valence of the vertex.
For
a picture see Figure 3.2. It suffices to show that
$${ \sum _ {x\in C\cap S_{{\hat {\Phi}}}}{\pm}f({{\hat \phi}}_{x})=0}\eqno
(3.2.1)
$$
for every interior vertex of $S_{{\hat {\Phi}}}$. Here
${\hat {\phi}}_x( S^1)={{\hat \Phi} (P\times \{x\})}.$

\p {\it Case 1}.  The valence of the interior vertex is one:
In this case it is easy to see that for $x\in S_{{\hat {\Phi}}}$,
near that vertex, the unique double point of ${\hat {\phi}}_{x}$ is at a kink.
So (3.2.1) is implied by the local integrability condition (1)

\p {\it Case 2}. The valence of the interior vertex is four: In this case the four poi
nts in
$C\cap S_{{\hat {\Phi}}}$ correspond to the four singular knots appearing in
the local
integrability condition (2) and one can show that (3.2.1) is
guaranteed by it. \qed

\medskip
\p {\bf Lemma 3.2.2.} {\sl Let $f: {\cal L}^{1} \longrightarrow {\cal R}$  be
a singular link invariant
and let  ${\Phi}: S^1 \longrightarrow  {{\cal M}^L(P, M)}$
a loop. Then $X_{\Phi}$ only depends on the free homotopy class
of  $\Phi$ in ${\cal M}^L(P, M)$.}
\medskip

\p {\bf Proof.} Let ${\Phi}^{'}$ be another closed homotopy in
almost general position such that $\Phi$,
$\Phi': S^1 \longrightarrow {\cal M}^L(P, M)$ are freely homotopic
loops in ${\cal M}^L(P, M)$. Then there exists
a homotopy ${{\Phi}_t: P \longrightarrow {\cal M}^L(P, M)}$ with
$t\in [0,1]$, such that $\Phi_0=\Phi$ and $\Phi_1=\Phi'$.

Let $\gamma$ be the path in $ {\cal M}^L(P, M)$ defined by
${\gamma}(t)= {\Phi}_t(L)$. After putting $\gamma$ in
almost general position we have
$$X_{\gamma\Phi'\gamma^{-1}}=X_{\gamma}+
X_{\Phi'}-X_{\gamma}=X_{\Phi'}.$$
Hence, we can assume that both $X_{\Phi}$ and $X_{{\Phi}^{'}}$ are
based at $L$ and the homotopy ${\Phi}_t$ is taken relatively $L$.
The homotopy ${\Phi}_t$ gives rise to a map
${\cal H}: P\times S^1 \times I \longrightarrow M$.
We cut the annulus $S^1\times I$ into a disc $D$ along a proper arc
$\alpha \subset S^1\times I$. Then, we have
$$X_{\partial D}={\pm}(X_{\Phi}-
X_{\Phi'}-X_{\alpha}+X_{\alpha}).$$
By Lemma 3.2.1 we obtain $X_{\partial D}=0$, and hence
$X_{\Phi}=X_{\Phi'}$. \qed

To continue, we first need to introduce some notation. Suppose that $P$ has $m$
components; that is $$P =\coprod_{i=1}^{m} P_i$$
\p where each $P_i$ is an oriented circle. Let $L: P
\longrightarrow M$ be a link.
Pick a base point $p_i \in P_i$ and
let $a_i$ denote the homotopy class of $L(P_i)$ in
$\pi_1 (M, L(p_i))$. Finally, we denote by $Z(a_i)$ the centralizer of
$a_i$ in $\pi_1(M, L(p_i))$.

\medskip
\p {\bf  Lemma 3.2.3.} {\sl Suppose that $M$ is a rational homology
3-sphere, with $\pi_2(M)=1$. Let $L: P \longrightarrow M$ be a link and let
$\Phi: P{\times}S^1\longrightarrow M$
be a closed homotopy from $L$ to itself. Moreover, assume that
either $Z(a_i)=\pi_1(M, L(p_i))$ or $Z(a_i)$ is finite for some $i=1,\dots, m$. Then
$X_{\Phi}=0$.}
\medskip

\p {\bf Proof.} Without loss of generality we may assume that
either $Z(a_1)=\pi_1(M, L (p_1))$ or $Z(a_1)$ is finite. Let
$L_1$ denote the restriction of $L$ on $P_1$ and let
${\Phi}_1$ denote the restriction of  ${\Phi}$ on $P_1\times S^1$.
We denote by ${\cal M}_1={\cal M}^{L_1}(P_1, M)$ the space of maps
$P_1 \longrightarrow M$, which are homotopic to $L_1$,
equipped with the compact-open topology. Let 
$$\pi= \pi_1 ({\cal M}^{L_1}(P_1, M),L_1).$$
One can see (see also Proposition 3.3 of [15]) that $\pi=Z(a_1)$.
Clearly, ${\Phi_1}$ represents an element in $\pi$.

Let $\Psi:P_1{\times}S^1\longrightarrow M$ be a loop in
${\cal M}_1$ based at $L_1$. We define ${\tilde {\Psi }}: P\times S^1
\longrightarrow M$ by
$$ {\tilde {\Psi }}|P_1 \times S^1 ={\Psi}$$
$$ {\tilde {\Psi }}(P^{'}\times S^1)=
\Phi (P^{'}\times S^1)$$
where  $P'=P\setminus P_1$. Then $\tilde\Psi$ is the closed homotopy from $L$
to itself.

After a small perturbation,
we can assume that there are only finitely many points
$x_1,x_2,\dots,x_n \in S^1$, ordered cyclicly according to the
orientation of $S^1$, so that ${\tilde {\psi}}_{x_i}\in {\cal L}^{1}$ and
${\tilde {\psi_x}}$ is equivalent to ${\tilde {\psi_y}}$ for all $x_i<x,y<x_{i+1
}$. (Here, ${\tilde {\psi}}_{t}$
denotes ${\tilde {\Psi}}(P\times \{t\})$).

Define
$${\chi}(\Psi):=X_{{\tilde {\Psi }}}=
\sum_{i=1}^{n}{\epsilon}_i f({{\tilde {\psi}}_{t_i})}\in {\cal R}$$
where ${\epsilon}={\pm 1}$ is determined as in (3.1.3), and
$\cal R$ is a ring which is torsion free as an abelian group.
Notice that $${\chi}({\Phi}_1)=X_{\Phi}$$

\medskip
\p {\it Claim:} The assignment $\Psi \longrightarrow {\chi}({\Psi})$
is a  group homomorphism
${\chi } : \pi \longrightarrow {\cal R}$.
\medskip

\p {\it Proof of the claim:} It is enough to show that ${\chi}({\Psi})$
is independent of the choice of the representative of $[\Psi]\in \pi$.
Let ${\Psi}_1 :  P_1{\times}S^1\longrightarrow M$ be another
loop in almost general position, which is homotopic to
$\Psi$. Then clearly ${\tilde {\Psi}}$ and ${\tilde {\Psi}}_1$
are homotopic loops in
${{\cal M}}^{L}(P, M)$, and the claim follows from Lemma 3.2.2.

By the assumption, we have either $\pi=\pi_1(M, L(p_1))$ or $\pi$ is finite, 
Since $\cal R$ is abelian, $\chi$ must factor through a finite abelian group (as $M$
is a rational homology 3-sphere).  Thus,  we must have $\chi =0$ since ${\cal R}$
is torsion free. In particular,
$${\chi}({\Phi}_1)=X_{\Phi}=0$$
as desired. \qed

\medskip
\p {\bf Corollary 3.2.4.} {\sl Assume that $M$ is as in Lemma 3.2.3 above. 
Let $L: P \longrightarrow M$ be a link and let $\Phi : P\times S^1 
\longrightarrow M$ be a closed homotopy
from $L$ to itself.

\p 1) If $L$ has a component
which is homotopically trivial in $M$, then $X_{\Phi}=0$.

\p 2) If $\pi_1(M)$ is finite, then $X_{\Phi}=0$.

\p 3) Assume that $M$ is a Seifert fibered space,
and that $L$ has a component which is homotopic to a regular fiber of 
the fibration. Then $X_{\Phi}=0$. }
\medskip

\p {\bf Proof.} 1)  and 2) follow immediately  from Lemma 3.2.3.

\p 3). By Lemma 32.8 of [8],  we know that the centralizer of a regular fiber
is $\pi_1(M)$. Hence the result follows from Lemma 3.2.3.\qed

Thus, we may assume from now on that $\pi_1(M)$ is infinite.

\medskip
\p {\bf \S3.3. Closed homotopy of links and essential tori}
\medskip

The purpose of this paragraph is study the topology of closed homotopy of links
thought of singular tori in 3-manifolds.  Since we are mainly interested in the
global integrability condition, which in general will be reduced to these special
cases discussed here, we may (and will) consider only 3-manifolds with infinite
$\pi_1$ by Corollary 3.2.4. 

Assume that $M$
is a Seifert fibered space with orbit surface $B$
and fiber projection $p: M\longrightarrow B$. Let
$\Phi: T=S^1\times S^1 \longrightarrow M$ be a closed homotopy of the knot 
$\Phi|S^1\times\{\ast\}$ which is 
vertical with respect to the given fibration. Let
$$\tilde M\longrightarrow M$$ 
be the covering space corresponding to
the cyclic normal subgroup generated by a regular fiber. We say that $\Phi
(S^1\times \{*\})$ {\it  doesn't wrap around the fibers of} $M$ if
$\Phi|S^1\times\{\ast\}$ lifts to $\tilde M$.

\p {\bf Lemma 3.3.1.} {\sl Let $M$ and $\Phi$ be as above. If $\Phi
(S^1\times \{*\})$ doesn't wrap around the fibers of $M$, then the closed
homotopy $\Phi$ is homotopic to another closed homotopy $\Phi'$ with the
following property: For every $x_1, x_2\in S^1$, there is a
homeomorphism $h^{12}:M\longrightarrow M$ such that

\p 1) $h^{12}=id$ outside of a regular neighborhood of $\Phi'(T)$ in $M$;

\p 2) $h^{12}(\phi_{x_1})=\phi_{x_2}$, where
$\phi_x=\Phi'|S^1\times\{x\}$;

\p 3) $h^{12}$ is isotopic to the identity map $id:M\longrightarrow M$.}
\medskip

See Lemma 3.9 in [12].

\medskip
\p {\bf Lemma 3.3.2. } {\sl Let $p:M\longrightarrow B$ be as above
and, in addition,  let $M$ be a rational homology 3-sphere. Let $\Phi :
T=S^1\times S^1 \longrightarrow M$ be an essential map.
Then, there exists a map $\Phi_1 : T \longrightarrow M$
which is homotopic to $\Phi$,
and a finite covering ${\tau}: S^1\times S^1 \longrightarrow S^1\times S^1$
such that the map $\Phi_1 \circ {\tau} : S^1\times S^1 \longrightarrow M$
can be extended to a map ${\hat {\Phi}}: S^1\times X
\longrightarrow M$. Here $X$ is a surface with
$\partial X=\{*\} \times S^1$.}
\medskip

\p {\bf Proof.}
 By Proposition 2.2.5, $\Phi $ is homotopic to a map
${\Phi}_1: T\longrightarrow M$ which is vertical with
respect to the fibration of $M$. Then, there exists a decomposition
$T=S^1\times S^1$ such that

\p a) ${\Phi}_1 (S^1\times \{*\})$ covers a regular fiber $h$, of $M$

\p b) We have $p({\Phi}_1 (\{*\}\times S^1))= p(T)$.

Let  $H$ (respectively, $Q$) denote the curve 
$S^1\times \{*\}$ (respectively, $\{*\}\times S^1$ on $T$,
and let $N$ be the cyclic normal subgroup  $\pi_1(M)$
generated by the regular fiber $h$.
Since $M$ is a rational homology 3-sphere, the abelianization of the
{\it Fuchsian } group $\Delta=\pi_1(M) / N$ is finite. Let $d$
be its order, and consider the $d$-fold covering
$\tau : {\hat T}: \longrightarrow T=H\times Q$,
corresponding to the subgroup $\Z \oplus d \Z$ of $\pi_1(T)=\Z \oplus \Z$.
Let ${\hat H}$ and ${\hat Q}$
denote the liftings, on ${\hat T}$, of $H$
and $Q$ respectively. Then, the map ${\hat Q}\longrightarrow B$
induced by $\Phi_1\circ \tau$ extends to a map $X \longrightarrow B$
for some compact surface $X$ with boundary  ${\hat Q}$.
This gives us a map $ \pi_1 (X) \longrightarrow \Delta$,
which in turn lifts to a map ${\hat {\Phi}}:
S^1\times X  \longrightarrow M$ (recall that $\pi_1(M)$
is an extension of $\Delta$ by $\Z\cong N$), with
${\hat {\Phi}}|S^1 \times {\partial X} = {\Phi_1\circ \tau}$.\qed

Let $P$ be a disjoint union of oriented circles and let $\Phi : P\times S^1 
\longrightarrow M$ be a
closed homotopy from a link $L: P \longrightarrow M$
to itself. Let $f:\cal L\longrightarrow\cal R$ be a singular link invariant
and let $X_{\Phi}$ be the quantity defined in (3.1.3). Suppose that
$P$ has
$m$ components, that is
$$P=\coprod_{i=1}^{m} P_i$$

\medskip
\p {\bf Lemma 3.3.3.} {\sl Assume that $M$ is a  Seifert fibered rational homology
3-sphere and let $P$, $\Phi$ be as above.
Moreover, assume that
$\Phi | P_i\times S^1$ is an essential map, for every $i=1,\dots, m$. Then there
exists a map ${\tilde {\Psi}}: P\times D^2 \longrightarrow M$ such that
$$X_{\partial\tilde{\Psi}}= a X_{{\Phi}}\eqno (3.3.1)$$
for some $a\in\Z$. Here, ${\partial\tilde{\Psi}}={\tilde {\Psi}}| P\times 
{\partial
D^2}$ and $D^2$ is a 2-disc. In particular, we have $X_\Phi=0$.}
\medskip

\p {\bf Proof.} Let $T_i=P_i\times S^1$
and let  ${\Phi}_i={\Phi}| T_i$, for $i=1,\dots, m$. Denote by $l_i$ (respectively,
$m_i$) the simple closed curve
$P_i\times \{*\}$ (respectively, $\{*\}\times S^1$) on $T_i$.

By Lemma 3.3.2, and after a homotopy to a vertical position, 
there exist  a finite covering ${\tau}_i : {\hat T_i}
\longrightarrow T_i$, such that 
${\Phi}_i\circ {\tau}_i$ extends
 to a map ${\hat {\Phi}_i}: S^1\times Y_i
\longrightarrow M$. Here $Y_i$ is a compact surface and
$ S^1\times {\partial Y_i}={\hat T_i}$.
Moreover all the ${\tau}_i$'s can be taken to be
of the same degree $d$. 

\medskip
\p {\it Case 1:} $d=1$ so that $\hat T_i=T_i$.  Notice, that this is always the case if
$H_1(M)=0$.

Recall, that the quantity
$X_{\Phi}$ doesn't change under homotopy (Lemma 3.2.2.),
and let $H_i$ (respectively, $Q_i$)
denote $S^1\times \{*\}$ (respectively, $\{*\}\times {\partial Y_i}$). 
Suppose that $l_i= a_i H_i +b_i Q_i$, for some $a_i$, $b_i \in \Z$.
We distinguish two subcases:

\medskip
\p {\it Subcase 1:} Suppose that $a_i\neq 0$ for every $i=1,\dots, m$.

Let $q_i : {\tilde T_i} \longrightarrow T_i$
be the covering of $T_i$ corresponding to the subgroup
$a_i {\Z} \oplus  {\Z}$ of $ \pi_1(T_i)=\Z\oplus \Z$. Let
${\tilde l_i}$,
${\tilde Q_i}$,${\tilde H_i}$ and ${\tilde m_i}$ denote the liftings of $l_i$, $
Q_i$, $H_i$ and $m_i$,
respectively. We have ${\tilde l_i}={\tilde H_i}+b_i {\tilde Q_i}$.

Each map $q_i$ extends to an $|a_{i}|$-fold covering,
$${\tilde q_i}: {\tilde l_i}\times {\tilde Y_i}
\longrightarrow S^1 \times Y_i$$
where ${\tilde Y_i}$ is a compact surface with $\partial {\tilde F_i}={\tilde Q_i}$,
and ${\tilde l_i}\times {\tilde Q_i}={\tilde T_i}$.

Let ${\tilde {\Phi}_i}={\hat \Phi_i}\circ{\tilde q_i}$ and let
$${\tilde {\Phi}}=\coprod_{i=1}^{m} {\tilde {\Phi}_i}$$

\medskip
\p {\it Claim:} We have that
$$X_{\partial\tilde {\Phi}}=a X_{ {\Phi}}$$
\p where $|a|=max\{|a_1|,\ldots,|a_m|\}$.
\medskip

\p {\it Proof of Claim:} Let ${q_i}_*$ denote the map induced by $q_i$ on the
fundamental groups. One can easily see that ${q_i}_*({\tilde m_i})=a_i m_i$
and $${\tilde Q_i}=c_i {\tilde l_i} + {\tilde m_i}\eqno (3.3.2)$$ for some $c_i
\in \Z$.
We identify the curves ${\tilde Q_i}$ by a common parameterization, and call the
result ${\tilde Q}$. The parameterization should be such that
corresponding points  on the ${\tilde Q_i}$'s
map, under the $q_i$'s, to the same point on the parameter space of $\Phi$.
By (3.3.2) this induces a common parameterization
of the curves ${\tilde m_i}$. Identify them and call the result ${\tilde m}$.
Now, ${\tilde {\Phi}}$ induces a map ${\tilde l}\times {\tilde m}
\longrightarrow M$, where
$${\tilde l}=\coprod_{i=1}^{m}{\tilde {l_i}}$$
We continue to denote this map by
${\tilde {\Phi}}$. Clearly, we have
$${\tilde {\Phi}}({\tilde l}\times \{x\})=
\coprod_{i=1}^{m} {\Phi_i}(P_i\times \{q_i(x)\})$$
for every $x$ on ${\tilde m}$. Notice
that each point, on the parameter space of ${\Phi}$, for which  ${ {\Phi}}(P)$
is not an embedding, corresponds to $|a|$ points
$x\in {\tilde m}$ for which
${\tilde {\Phi}}({\tilde l}\times \{x\})$
is not an embedding. Now, the claim follows easily.
Let us finally observe that, because of (3.3.2),
the quantity $X_{\partial\tilde {\Phi}}$ doesn't change
if we replace the parameter space ${\tilde m}$, by ${\tilde Q}$.

To continue with the proof of the lemma, we choose a collection of proper arcs
$${\{{{\alpha}_i^j}\}}_{j=1}^{m_i} \subset {\tilde Y_i}$$
such that: a) each $ {\tilde Y_i}$ if cut along the ${\{{{\alpha}_i^j}\}}$'s
becomes a disc (as each $\tilde Y_i$ can be chosen as connected),
and b) the end points of the
${\{{{\alpha}_i^j}\}}$'s avoid the points for which ${\tilde {\Phi}}(\tilde l\times
\{*\})$ is not an embedding. Let $\Gamma_i$ denote
the space obtained by cutting ${\tilde l_i}\times {\tilde Y_i}$
along the collection of annuli
$${\{{{A}_i^j}\}}_{j=1}^{m_i}$$
where $A_i^j={\tilde l_i}\times {{{{\alpha}_i^j}}}$. Let us
denote by
${\tilde \Psi}_i$ the map induced on $\Gamma_i$ ,
by ${\hat \Phi_i}\circ{\tilde q_i}$. Finally, let
$$\Gamma=\coprod_{i=1}^{m} {\Gamma_i}$$
and let
$${\tilde \Psi}=\coprod_{i=1}^{m}{\tilde \Psi}_i$$
The map induced on $\Gamma$ by
${\tilde \Psi}$, is the desired map.

\medskip
\p {\it Subcase 2:} Suppose that $a_i=0$ for some $i=1,\dots,m$. 

Suppose for example that $a_1=0$. Then $\Phi(l_1)$ doesn't wrap around the
fibers of $M$.  By Lemma 3.3.1,  we may assume that 
$\Phi_1 (P_1\times \{x_1\})$ and
$\Phi_1 (P_1\times \{x_2\})$ are isotopic
for every $x_1$ and $x_2 \in S^1$. Let
$$P'=\coprod_{j\neq 1} P_j$$
\p and let
$${\Phi}'=\coprod_{j\neq 1} \Phi_j$$

Observe that  $\Phi (P\times \{*\})$ is not an embedding if either
$ {\Phi}'(P'\times \{*\})$ is not an embedding,
or $\Phi_1 (P_1\times \{*\})$
intersects with some $\Phi_j (P_j\times \{*\})$.
We may change $\Phi'$ by composing it with the inverse of the isotopy of
$\Phi_1(P_1\times\{*\})$. Hence,
$X_{\Phi}$ doesn't change if we  assume that $$\Phi_1 (P_1\times \{x\})=\Phi_1
(P_1\times
\{x_0\})$$
\p where $x_0$ is fixed and $x$ runs on $\{*\} \times S^1$. Hence,
${\Phi_1}$ extends to a map ${\hat {\Phi}_1}: S^1\times D_1
\longrightarrow M$, where $D_1$ is a disc.
Then we proceed as in Subcase 1 above. 

\medskip
\p {\it Case 2:} Now assume that $d>1$,
where $d$ is the common degree of the coverings 
${\tau}_i: {\hat T_i}
\longrightarrow T_i$.

In this case
we apply Lemma 3.3.2 to a suitable covering of $T_i$ rather that $T_i$
itself. More precisely suppose that $T_i$ is in vertical position, and
that $l_i= a_i H_i +b_i Q_i$, for some $a_i$, $b_i \in \Z$.

If $a_i\neq 0$ for every $i=1,\dots, m$,
let
$$q_i : {\tilde T_i} \longrightarrow T_i$$
 be the covering of $T_i$ corresponding to the subgroup
$a_i {\Z} \oplus  {\Z}$ of $ \pi_1(T_i)=\Z\oplus \Z$. Let
${\tilde l_i}$,
${\tilde Q_i}$,${\tilde H_i}$ and ${\tilde m_i}$ denote the liftings of $l_i$, $
Q_i$, $H_i$ and $m_i$,
respectively. We have ${\tilde l_i}={\tilde H_i}+b_i {\tilde Q_i}$
and we may choose ${\tilde L_i}$ and ${\tilde Q_i}$
as a system of generators of $\pi_1 (\tilde T_i)$.
Let
$${\tilde {\Phi}}=\coprod_{i=1}^{m} { {\Phi}_i}\circ q_i$$

In view of the claim above, it is enough to prove the assertion in the statement
of the lemma for ${\tilde {\Phi}}$. Now let
 ${\tilde \tau}_i: {\tilde T_i}^*\longrightarrow {\tilde T_i}$
 be coverings as in the proof of Lemma 3.3.2, and let
 $d$ be their common degree. 
 Let ${\Phi}_i^*={ {\Phi}_i}\circ q_i \circ {\tilde \tau}_i$
 and let 
 $${{\Phi}}^*=\coprod_{i=1}^{m} { {\Phi}_i^*}$$
 One can see that
$$X_{\Phi^*}= d X_{\tilde \Phi}$$
and proceed as in Case 1 above to prove the desired assertion for ${\Phi^*}$.

If $a_i=0$ for some $i=1,\ldots,m$, then we proceed
as in  Subcase 2 above. \qed

\medskip
\p {\bf Lemma 3.3.4.} {\sl  Assume
that $M$ is an oriented 3-manifold with $\pi_2(M)=1$.
Let $\Phi$ be a closed homotopy such that $$\Phi | P_i\times S^1$$
is an inessential map, for every $i=1,\dots, m$. Then $X_\Phi=0$.}
\medskip

\p {\bf Proof.} As before, we denote  $T_i=P_i\times S^1$, ${\Phi}_i={\Phi}| T_i$,
$l_i=P_i\times \{*\}$, and $m_i=\{*\}\times S^1$ for $i=1,\dots,m$.

We have assumed that $\pi_1(M)$ is infinite. So it is torsion free (Theorem
9.8 of [6]).  Hence, ${\Phi}_i$ extends to a map ${\hat {\Phi}_i}: S^1\times D_i
\longrightarrow M$ , where $D_i$ is a 2-disc and
$ S^1\times {\partial D_i}=T_i$. Let $H_i$ (respectively, $Q_i$)
denote $S^1\times \{*\}$ (respectively, $\{*\}\times {\partial D_i}$). Suppose that
$l_i= a_i H_i +b_i Q_i$, for some $a_i$,
$b_i \in \Z$.

If $a_i=0$ for some $i=1,\ldots,m$ then one
of the components of the link $\Phi{(P\times \{*\})}$
is homotopically trivial and the conclusion follows from Lemma 3.2.3.
In the case  that $a_i\neq 0$ for every $i=1,\dots, m$, one proceeds as in Case 1
of
 the proof of Lemma 3.3.3 to obtain a map ${\tilde {\Psi}}:
P\times D^2 \longrightarrow M$ such that
$$X_{\partial{\tilde{\Psi}}}=a X_{{\Phi}}\eqno (3.3.3)$$
for some  $a\in \Z$. The desired conclusion then follows immediately. \qed

\medskip
\p {\bf \S3.4. Completing the proof of Theorem 3.1.2}
\medskip

Recall that $M$ is a rational homology 3-sphere which is either atoroidal or a 
Seifert fibered space. Also, $\pi_1(M)$ is infinite.
As before, $P$ is a disjoint union of oriented circles and 
$$\Phi : P\times S^1 \longrightarrow M$$
is a
closed homotopy from some link $L: P \longrightarrow M$
to itself. Let $f: {\cal L}^{1} \longrightarrow {\cal R}$ be a singular link 
invariant as in the statement of theorem 3.1.2.
We have to show that
 $$X_{\Phi}=0 \eqno (3.4.1)$$
where $X_{\Phi}$ is the signed sum of
values of $f$ around $\Phi$ defined
in (3.1.3).

If $M$ is atoroidal then by Lemma 3.3.4,  (3.4.1) holds.

Suppose that $M$ is a  Seifert fibered space.  Let $E$  (resp. $I$) denote
the set of components of $P\times S^1$ on which
$\Phi$ is  essential (resp.  inessential). If $E$ or $I$ is
empty the claim follows from Lemmas 3.3.3 or 3.3.4.  

We first notice that as in the proof of Lemma 3.3.4, 
if there is  a component in
$I$ which is a closed homotopy of a homotopically trivial knot, then $X_\Phi=0$
by Corollary 3.2.4.  Otherwise, we have the following claim.

\medskip
\p {\it Calim:} There exists a map ${\tilde {\Psi}}: P\times D^2 \longrightarrow M$
such that
$$X_{\partial{\tilde{\Psi}}}=a X_{{\Phi}}$$
for some  $a\in \Z$. In particular, we $X_\Phi=0$.

\medskip
\p {\it Proof of the claim:} The proof is very similar to the
proofs of Lemmas 3.3.3 and 3.3.4. The only difference is that,
for example, when cutting each surface $Y_i$ to discs (see the proof of Lemma
3.3.3),  we have to make sure that the end points of cutting arcs 
also avoid singular links where the double points are other than intersections
between components of $L$ in $E$. The details are left to the reader.

This finishes the proof of Theorem 3.1.2. \qed

\medskip
\p {\bf Remark 3.4.1.} A more general statement, than that
of Theorem 3.1.2, was proved in [15] for links in manifolds
with finite $\pi_1$ and in [12] for knots in irreducible 3-manifolds,
as a first step in understanding {\it finite type} invariants
for knots in these manifolds.
\bigskip

\p{{\bf \S 4. An intrinsic definition of the HOMFLY power series }}
\medskip
\p {\bf \S4.1. Preliminaries}
\medskip

It is known ([11]) that the 2-variable Jones (or HOMFLY) polynomial ([5,10]) 
for links in $\R^3$ or $S^3$ is equivalent to sequence 
${\{J_n=J_n(t)\}}_{n\in \Z}$ of 1-variable
Laurent polynomials. They are completely determined
 by the following skein relations:
$$J_n(U)=1\eqno (4.1.1)$$
$$t^{{n+1}\over 2}J_n({L_{+}})
-t^{-{{n+1}\over 2}} J_n({L_{-}})=(t^{1\over 2}-t^{-{1\over 2}})  J_n(L_{o})
\eqno (4.1.2)$$
 where $L_{+}$, $L_{-}$, $L_o$ are  the resolutions
of a singular link $L_{\times}\in {\cal L}^{(1)}$ described in 3.1.
In this context, the original Jones polynomial is $J_{-3}$. 

Notice that the initial value $J_n(U)=1$ is not essential. Any choice of the initial value
together with (4.1.2) will determine a unique $J_n$.

 Let $$u_n(t)= {{{t^{{n+1}\over 2}-t^{-{{n+1}\over
2}}}\over {t^{1\over 2}-t^ {-{1\over 2}}}}}.$$
By (4.1.2) one obtains
$$ J_n(L\coprod U) = u_n(t)\  J_n(L) \eqno (4.1.3)$$
where the link $L\coprod U$ is obtained from $L$ by adding an unknotted and
unlinked component $U$.

The coefficients of the power series $J_n(x)$, obtained from $J_n(t)$
by substituting $t=e^ x$, are invariants of {\it finite type} ([2,3]).
In the theorem below we reverse this procedure, and guided by (4.1.2) we will 
construct inductively power series invariants for links in 3-manifolds generalizing the
$J_n(x)$'s.
 
\medskip
\p {\bf \S 4.2. The construction of the invariants}
\medskip

Assume that $M$ is an orientable, rational homology
3-sphere which is either atoroidal
or a Seifert fibered space.
For every $n\in \Z$, we will construct a sequence of knot invariants
$$v_n^0,\ v_n^1, \dots, \ v_n^m,\dots$$
such that the formal power series
$$ J_{\{ M, n\}}(L)= \sum_{m=0}^{\infty}v_n^m(L)x^m$$
satisfy (4.1.2), under the change of variable $t={e}^x$, for every
$L \in {\cal L}$.

We will construct our invariants inductively
(induction on $m$) by using Theorem 3.1.2. More precisely, each $v_n^m$ 
is going
to be obtained by integrating a suitable singular link invariant
determined by the $v_n^j$'s with $j<m$.

Recall that a link
invariant
obtained by integrating a singular link invariant
is well defined up to
a collection of ``integral constants" (see the beginning of
the proof of Theorem 3.1.2).
This means that in order to define
$v_n^0, v_n^1, \dots, v_n^m,\dots$ uniquely, we need
to make a choice of ``initial links".

Let $L$ be an $n$-component link and recall from \S3 that
${\cal M}^L$ denotes the space of maps
\p $\coprod S^1\longrightarrow M$
which are  homotopic $L$. The spaces ${\cal M}^L$
corresponding to links with $n$ components are in
one to one correspondence with
the  unordered $n$-tuples
of conjugancy classes in $ \pi=\pi_1(M)$. In every such space we will fix, 
once and for all, a link $CL$ and call it a
{\it trivial link}.
If $CL$ has $k$ components  which are homotopically
trivial, our choice will be such that $CL=CL^*\coprod U^k$, where $U^k$
is the standard unlink with $k$ components in a small ball neighborhood
disjoint from $CL^*$. Let $\cal C\cal L$ be the set of trivial links and
$\cal C\cal L^*$ be the set of trivial links with homotopically
non-trivial components.

Notice that when $M$ is simply connected there is a natural
choice of trivial links. Namely, one chooses each $CL$ to be
the unlink $U^k$. Thus $\cal C\cal L^*=\emptyset$ in this case.

\medskip
\p {\bf Theorem 4.2.1.} {\sl Let $M$, ${\cal L}$ and 
${\cal C\cal L^*}$ be as above.
There exists a  {\it unique} sequence of complex valued 
link invariants $v_n^0, v_n^1, \dots, v_n^m,\dots$, 
with given values on the links in
${\cal C\cal L^*}\cup\{U\}$,
such that if we define a formal power series
$$ J_{\{ M, n\}}(L)= \sum_{m=0}^{\infty}v_n^m(L)x^m$$
for $L \in {\cal L}$ then
$$t^{{{n+1}\over 2}}J_{\{ M, n\}}({L_{+}})
-t^{-{{{n+1}}\over 2}} J_{\{ M, n\}}({L_{-}})=
(t^{1\over 2}-t^{-{1\over 2}})  J_{\{M,n\}}(L_{o})\eqno (4.2.1)$$
where $t={e}^x=1+x+{{x^2}\over 2}+\dots$}
\medskip

\p {\it Notation:} To simplify our notation, and
throughout this proof, we will write $J_n$ instead of
$ J_{\{ M, n\}}$.
\medskip

\p {\bf Proof.} By our assumption, the values
$v_n^0(CL^*), v_n^1(CL^*), \dots, v_n^m(CL^*), \dots$ are given
for every $CL^*\in {\cal C\cal L}^*$. Hence, we can form the power series
$J_n(CL^*)$. Also, we may form $J_n(U)$ using the given values $v_n^m(U)$'s.

Guided by (4.1.3) we define
$$ J_{n}(CL\coprod U) = u_n(t)\  J_n(CL)\eqno(4.2.2)$$
where
$$u_n(t)={t^{{n+1}\over 2}-t^{-{{n+1}\over 2}}}\over{t^{1\over 2}-t^{-{1
\over 2}}}.\eqno (4.2.3)$$
Thus, $J_n$ has been defined on all trivial links.

We define the link invariant $v_n^0$ by
$$v_n^0(L)=v_n^0(CL),$$
where
$CL$ is the trivial link homotopic to $L$.

Inductively, suppose that the invariants
$v_n^0, v_n^1, \dots, v_{n}^{m-1}$
have been defined such that if we let
$$J_{ n}^{(m-1)}(L)= \sum_{i=1}^{m-1}v_n^i(L)x^i,$$
then
$$ J_{n}(L\coprod U) = u_n(t)\  J_n(L) \ \ mod\ \ x^m, \eqno (4.2.4)$$
$$t^{{n+1}\over 2}J_n ^{(m-1)}({L_{+}})
-t^{-{{{n+1}}\over 2}} J_n^{(m-1)}({L_{-}})=(t^{1\over 2}-t^{-{1\over 2}})  J_n^
{(m-1)}(L_{o})\ \  mod\ \ x^m. \eqno (4.2.5)$$
>From (4.2.5) we obtain
$${\eqalign {& J_n^{(m-1)}(L_{+})-J_n^{(m-1)}({L_{-}})=\cr
 &\ \ (t^{-(n+1)}-1)J_n^{(m-1)}({L_{-}})+t^{-{{n+1}\over 2}}
(t^{1\over 2}-t^{-{1\over 2}})J_n^{(m-1)}(L_{o})\ \  mod\ \ x^m \cr }}$$
which leads us to define
$${\eqalign {& J_n^{(m)}({L_{\times}}):=\cr
&\ \ (t^{-(n+1)}-1)J_n^{(m-1)}({L_{-}})+t^{-{{n+1}\over 2}}
(t^{1\over 2}-t^{-{1\over 2}})J_{n}^{(m-1)}(L_{o})\ \  mod\ \  x^{m+1}.\cr }
}\eqno (4.2.6) $$

This is a polynomial of degree $m$, with trivial constant coefficient.
The coefficients of $x^i$, $i=1$, $2$, $\ldots$,$m-1$, in this polynomial
are singular link invariants derived from $v_n^i$, $i=1$, $2$, $\ldots$,$m-1$.
However, the coefficient of $x^m$ is a {\it new} singular link invariant.
We are going to prove that it is derived from a knot invariant
by using Theorem 3.1.2. For that we need to check  that the local
integrability conditions (1) and (2) are satisfied. It is enough to check them
modulo $x^{m+1}$. In what follows the symbol $\equiv$ will denote calculation
modulo $x^{m+1}$.

To check (1), suppose we start with a singular link  $L^1\in {\cal L}^1$
that contains a kink {{\epsfysize=0.075in\epsfbox{F-4.ps} }}.
Let $L$ be the link obtained by resolving the double
point of $L^1$ and let $U$ be {\it unknot} in $M$. Then we have
$${\eqalign {& J_n^{(m)}(L^1)  \cr
&\equiv (t^{-(n+1)}-1)J_n^{(m-1)}(L)+
t^{-{{{n+1}}\over 2}}
(t^{1\over 2}-t^{-{1\over 2}}) \ J_n^{(m-1)}(L\coprod U )\cr
& \equiv \big[(t^{-(n+1)}-1)+t^{-{{{n+1}}\over 2}}(t^{1\over 2}-t^{-{1\over 2}})
\  u_n(t)\big]J_n^{(m-1)}(L) \cr
&\equiv 0\cr}}$$

To check (2), we calculate
$$\eqalign { & J_n^{(m)}(L_{{\times }+})-
J_n^{(m)}(L_{{\times }-})\cr
& \equiv (t^{-(n+1)}-1)J_n^{(m-1)}(L_{-+})+
t^{-{{{n+1}}\over 2}}
(t^{1\over 2}-t^{-{1\over 2}})J_n^{(m-1)}(L_{o +})\cr
&- (t^{-(n+1)}-1)J_n^{(m-1)}(L_{- -})-
t^{-{{{n+1}}\over 2}}
(t^{1\over 2}-t^{-{1\over 2}})J_n^{(m-1)}(L_{o -}) \cr
&\equiv (t^{-(n+1)}-1)
\big[ (t^{-(n+1)}-1)J_n^{(m-1)}(L_{- -})+
t^{-{{{n+1}}\over 2}}
(t^{1\over 2}-t^{-{1\over 2}})J_n^{(m-1)}(L_{- o}) +o(x^m)\big]\cr
&+ t^{-{{{n+1}}\over 2}}
(t^{1\over 2}-t^{-{1\over 2}})
\big[ (t^{-(n+1)}-1)J_n^{(m-1)}(L_{o -})+
t^{-{{{n+1}}\over 2}}
(t^{1\over 2}-t^{-{1\over 2}})
J_n^{(m-1)}(L_{o o}) +o(x^m)\big] \cr
&\equiv {(t^{-(n+1)}-1)}^2 J_n^{(m-1)}(L_{- -}) +{\big[t^{-{{{n+1}}\over 2}}
(t^{1\over 2}-t^{-{1\over 2}})\big]}^2 J_n^{(m-1)}(L_{o o})\cr
&+ (t^{-(n+1)}-1) t^{-{{{n+1}}\over 2}}
(t^{1\over 2}-t^{-{1\over 2}})\big[ J_n^{(m-1)}(L_{o -})+
J_n^{(m-1)}(L_{- o}) \big].\cr
 }$$

Since the result is symmetric with respect to the two double points
we deduce that
$$J_n^{(m)}(L_{{\times}+})-J_n^{(m)}(L_{{\times}-})\equiv
J_n^{(m)}(L_{+{\times}})-J_n^{(m)}(L_{-{\times}})$$

Thus, the singular link invariant defined in (4.2.6)
is induced by a link invariant. Using the given
values $\{v_n^m(CL):\ CL\in {\cal CL}\}$,
we can define a link invariant $v_n^m$, such that if we let
$$J_{ n}^{(m)}(L)= \sum_{i=1}^{m}v_n^m(L)x^i$$
we have
$$J_n^{(m)}({L_{+}})-J_n^{(m)}({L_{-}})=J_n^{(m)}({L_{\times}})$$
for $L\in {\cal L}$ and $L_{\times }\in {{\cal L}^1}$.
Therefore the invariant $J_n^{(m)}$ defined in this way satisfies the 
inductive hypothesis (4.2.5).

Now, a straightforward calculation shows that
$$J_n^{(m)}(L_{{\times}}\coprod U)\equiv u_n(t)\  J_n^{(m)}(L_{{\times}})$$
which together with (4.2.2)
shows that $J_n^{(m)}$ satisfies the inductive hypothesis (4.2.4).

To finish our proof we need to show uniqueness. Inductively, we assume that
$v_n^0, v_n^1, \dots, v_{n}^{m-1}$
are uniquely determined by (4.2.1)
and their values on ${\cal CL}$, for every $n\in \Z$.
Then, the conclusion for $v_n^{m}$ follows from the fact that
$$v_n^m(L)=v_n^m(CL)+ \sum_{i=1}^s {\pm} v_n^m (L_i)$$
where $L_1, \dots, L_s$ are singular links in ${\cal L}^1$,
and $CL$ is the representative of $L$ in ${\cal CL}$. \qed

To illustrate how the power series $J_n=J_{\{ M, n\}}$
depend on the choice of ${\cal CL}$,
let us restrict ourselves to the case of knots. Let
$\cal C$ denote the set of homotopy classes of loops in $M$.
Then ${\cal CK}$ will contain exactly one knot, say $K_C$
from every $C\in {\cal C}$. Moreover, we suppose that
the initial values of $J_n$ are all equal to $1$.
To stress out the dependence on $K_C$, let 
$$J_n^{K_C}(K):=J
_n(K)$$
for every $K\in {\cal M}^{K_c}(S^1,M)$.

Let
$ {\{ \tilde K_C\}}_{C\in {\cal C}}$ be a different choice of 
trivial knots.
Then one can see that $J_n(K)$ is well defined up to a multiplicative 
constant in the ring
of formal power series with complex coefficients. More precisely,

\medskip
\p {\bf Proposition 4.2.2.} {\sl We have
$$J_n^{\tilde K_C}(K)=J_{n}^{\tilde K_C}(K_C)J_{n}^{ K_C}(K)$$
for every  $K\in {\cal M}^{K_c}(S^1,M)$ and $c\in {\cal C}$.}
\medskip

\p {\bf Proof.} Follows immediately from the definition. \qed

\medskip
\p {\bf Remark 4.2.4.}
Recall that an invariant $f$, is called of {\it finite type}
if there exists an integer $m$, such that the singular
link invariant derived from $f$ is zero on all singular links with more than
$m$ double points. One can see that the invariants
$v_n^0,\ v_n^1, \ldots \ v_n^m,\ldots$ constructed above
are of {\it finite type}. It would be interesting to find a direct relation
of the invariants constructed here with these coming from the 
$SU(N)$-perturbative
Chern-Simons theory ([1,19]).
\bigskip

\p{\bf  \S 5. The convergence of the power series}
\medskip

The purpose of this section is to show that, with an appropriate choice of the inital
values, we have 
$$J_{\{M,n\}}(L)=J_{\{S^3,n\}}(L')$$
for a certain link $L'$ in $S^3$. See Theorem 5.2.3. Since $J_{\{S^3,n\}}$ is a Laurent
polynomial in 
$$v=t^{-{{n+1}\over 2}}\ \ \ \ \ \hbox{and}\ \ \ \ \ z=t^{1\over 2}-t^{-{1\over 2}},$$
so does $J_{\{M,n\}}$.
Hence, we obtain an actual generalization
of the 2-variable Jones-polynomial for links in atoroidal and Seifert fibered rational
homology 3-spheres.

\medskip
\p {\bf 5.1. Some special cases}
\medskip

The reader might have noticed (see also [15]) that for $M=S^3$, Theorem 4.2.1
gives an intrinsically 3-dimensional proof of the existence
of the HOMFLY ``series". To see that it converges to a Laurent polynomial in $v$ and
$z$, though, we have to use some special features of $S^3$. Namely, every link in
$S^3$ has a plane projection and when (4.2.1) can be applied to simplify that plane
projection. So, if the initial value, $J_{\{S^3,n\}}(U)$ in this case, is chosen to be a Laurent
polynomial in $v$ and $z$, $J_{\{S^3,n\}}(L)$ will be a Laurent polynomial for every $L$. 
Thus, the main question we are facing now is the following.

\medskip
\p {\bf Question 5.1.1.} {\sl Do the power series $J_{\{M,n\}}(t)$
converge to Laurent polynomials or rational functions
in other manifolds M, besides $S^3$?}
\medskip

Theorem 5.2.3 gives a partial answer to this question. But we feel that it is still not 
clear how the convergence problem is related with the topology of the underlying
manifold. 

Let us see a special case first. Let us define
$${\nabla}_L={\nabla}_L(t):=J_{\{M,-1\}}(L)(t^2)$$
By (4.2.1) we see that ${\nabla}_L$
satisfies the following skein relation
$$\nabla_{L_+}-\nabla_{L_-}=(t-t^{-1})\nabla_{L_o}\eqno (5.1.1)$$
By (5.1.1) and the results in [4],  we see that, possibly after some normalization,
${\nabla}_L$ is the {\it Conway potential function}, which
is known to be a rational function in $t$. More precisely,
it is known that
$${\nabla}_L(t)={{\Delta_L(t^2)}\over {t-t^{-1}}}\eqno (5.1.2)$$
if $L$ is a knot, and
$${\nabla}_L(t)={\Delta_L(t^2,\ldots, t^2)}\eqno(5.1.3)$$
if $L$ is a link. Here $\Delta_L$ denotes the Alexander polynomial of $L$,
and the number of variable entries in the right hand side of (5.1.3)
is the same as the number of components of $L$. For more details see [4].

\medskip
\p {\bf 5.2. The existence of the 2-variable Jones polynomial}
\medskip

In this paragraph we will use the fact that $M$ is obtained
by surgering $S^3$ along a link ([14]) together with the discussion in
the beginning of 5.1 to show that the $J_{\{M, n\}}(t)$'s
can be suitably normalized so that they converge to Laurent
polynomials in $v$ and $z$, giving thereby a positive answer to Question 5.1.1.
To simplify our notation we will write
$J_n$ instead of $J_{\{ S^3, n\}}$.
We need the following
\medskip
\smallskip
\p {\bf Lemma 5.2.1.} {\sl There exists links $L^1\subset M$ and
$L^0\subset S^3$, such that $M\setminus L^1$ is homeomorphic
to $S^3\setminus L^0$. We will write $M\setminus L^1\cong S^3\setminus L^0$.}
\medskip
\p {\bf Proof.} There exists a framed link $L\subset S^3$ such that
 the result, say ${\chi(L)}$, of surgery  of $S^3$ along $L$
is homeomorphic to $M$. Hence,
there exists a homeomorphism $h_0: {\chi(L)} \longrightarrow M$.
Then $L^0=L$ and $L^1=h_0(L)$ are as desired.\qed
\medskip
Let us fix $L^1$ and $L^0$ as in the lemma above  and let
$M_1=M\setminus L^1$. Let $L\subset M$ be a link.
After isotopy we may assume that $L\cap L^1 =\emptyset$ and hence
$L\subset M^1$. Moreover, we have the following lemma whose proof is quite obvious.

\medskip
\p {\bf Lemma 5.2.2.} {\sl Let $L_1$, $L_2\subset M^1$ be two links
which are homotopic to each other in $M$. Then there exist a link ${\tilde L}$ in $M_1$
such that

1) ${\tilde L}$ is isotopic to $L_2$ in $M$; and

2) ${\tilde L}$ is homotopic to $L_1$ in $M_1$.}
\medskip

Now, let us fix a
homeomorphism $h: M_1\longrightarrow S^3\setminus L^0$.  Suppose that we
have chosen every
$CL^*$ to be in $M^1$. By Theorem 4.2.1, there exists a  {\it unique} sequence of
numerical link invariants
$v_{\{M,n\}}^0, v_{\{M,n\}}^1,
\dots,  v_{\{M,n\}}^m,\dots$ such that the formal power series
$$ J_{\{ M, n\}}(L)= \sum_{m=0}^{\infty}v_{\{M,n\}}^m(L)x^m$$
has the following properties:

1) $ J_{\{ M, n\}}(CL^*)=J_n(h(CL^*))$ for every $CL^*\in {\cal CL}^*$; 

2) $J_{\{M,n\}}(U)=1$; and

3)  $J_{\{M,n\}}$ satisfies the skein relation (4.2.1).

Here $J_n=J_n(t)=\sum_{m=1}^\infty v_n^mx^m$ is the power series obtained from
the polynomial
$J_n(t)$ by the substituting $t=e^x$. We have,

\medskip
\p {\bf Theorem 5.2.3.} {\sl For every
link $L\in {\cal L}(M)$ and $n\in \Z$, we have
 $$ J_{\{ M, n\}}(L)= J_n(h(L))$$
Thus, since the power series $J_n$ converges to a Laurent polynomial in
$$v=t^{-{{n+1}\over 2}}\ \ \ \ \ \hbox{and}\ \ \ \ \ z=t^{1\over 2}-t^{-{1\over 2}},$$
so does $J_{\{M,n\}}$.}
\medskip

\p {\bf Proof.}  Let us assume inductively that the invariants
$v_{\{M,n\}}^0, v_{\{M,n\}}^1, \dots, \ v_{\{M,n\}}^{m-1} $ satisfy
$$v_{\{M,n\}}^i(L)=v_n^i(h(L))\eqno (5.2.1)$$
for every  $L\in M^1\subset M$  and $i=1, \dots, m-1$.
We will show that (5.2.1) is true for $i=m$.

By the assumption, we have 
$$J_{\{M,n\}}(CL)=J_n(h(CL))$$
for every $CL\in{\cal CL}$. Therefore, 
 $$v_{\{M,n\}}^m(CL)=v_n^m(h(CL))\eqno (5.2.2)$$
for every $CL\in {\cal CL}$.
Now, let $L\subset M^1\subset M$ be a link. By our assumption.
$L$ is homotopic to some $CL\in {\cal CL}$. By Lemma 5.2.2
we may assume that $L$ is homotopic to $CL$ in $M^1$.
Hence (see the proof of Theorem 4.2.1)
there exist singular links  $L^1$,$\ldots$, $L^s$
in $M^1$ such that
$$v_{\{M,n\}}^m(L)=v_{\{M,n\}}^m(CL)+ \sum_{i=1}^s {\pm} v_{\{M,n\}}^m(L^i)\eqno
 (5.2.3)$$
And we also have
$$v_n^m(h(L))=v_n^m(h(CL))+ \sum_{i=1}^s {\pm} v_n^m(h(L^i))\eqno (5.2.4)$$

Let $L_-^i$ (respectively, $L_o^i$) denote the link obtained from $L^i$
by resolving its double point negatively
(respectively surgering the double point in a way consistent
with the orientation). They are links in $M^1$. By (4.2.6) we see that $v_{\{M,n\}}^m(L^i)$
is a linear combination of $v_{\{M,n\}}^j(L^i_-)$'s
and $v_{\{M,n\}}^j(L^i_o)$'s  with $j=0, \dots, m-1$.
Hence, by the induction hypothesis we have,
$$v_{\{M,n\}}^j(L^i_-)=v_{n}^j(h(L^i_-))$$
and
$$v_{\{M,n\}}^j(L^i_o)=v_{n}^j(h(L^i_o))$$
for every $j=0, \dots, m-1$. One can see that
$h(L^i)$, $h(L^i_-)$ and $h(L^i_o)$ are related (in $S^3$) in
the same way that $L^i$, $L_-^i$ and $L_o^i$ are related in $M^1$.
Hence, we obtain
$$ v_{\{M,n\}}^m(L^i)= v_{n}^m(h(L^i))\eqno (5.2.5)$$
Now, the desired conclusion follows immediately from
(5.2.2) to (5.2.5).\qed

\medskip
\p {\bf Proof of Theorems A and B.} 
Combine Theorem 4.2.1 and Theorem 5.2.3, to obtain Theorem A. 
Theores B will follow immediately from
Theorem A since ${\cal CL}^*=\emptyset$ when $\pi_1(M)=1$.\qed
\bigskip

%\magnification=1200
{\centerline{\bf References}}
\medskip
%\item{[AS]} S. Axelrod and I.M. Singer, { Chern-Simons
%perturbation theory}, M.I.T. preprint.

%\item {1.}  {J.S Birman}, {\it New  points of  view
 % in knot theory}, Bulletin of AMS, vol. {\bf 28}, no {\bf 2} (1993)
%253--287.

%\item{[B--K]} J.S. Birman and T. Konenobou, {\sl
%Jones' braid-plat formula and a new surgery triple},
%Proc. of AMS {\bf 102}, 1988.

\item{1.} D. Bar-Natan, {\it Perturbative aspects of the
Chern-Simons topological quantum field theory}, Ph.d thesis,
Princeton University, 1991.

\item{2.} { {D. Bar-Natan}},
{\it On the Vassiliev knot invariants},
Topology, {\bf 34}(1995), 423-472.

\item {3.} { {J.S. Birman} { and}  {X.S. Lin}},
{\it  Knot polynomials and Vassiliev invariants}, Invent. Math.,
{\bf 111}(1993),225--270.

%\item{[BW]} J. S. Birman and H. Wenzl, {\sl Braids, link
%polynomials and a new algebra}, Trans. AMS {\bf 313} (1989)
%249-273.

\item{4.} S. Boyer and D. Lines,
{\it A Conway potential function for Rational
Homology Spheres}, Proc. Edinburgh Math. Soc. (2), {\bf 35}(1992), 53--69.

\item{5.} P. Freyd, J. Hoste, W.B.R. Lickorish,
K. Millet, A. Ocneanu and D. Yetter, {\it A new
polynomial invariant of knots and links}, Bull. AMS {\bf
12}(1985), 239-246.

\item{6.} J. Hempel, {\it 3--manifolds}, Ann. of Math.
Studies, vol.  {\bf 86}, Princeton University Press, 1976.

\item{7.} J. Hoste anf J. Przytyski, {\it A survey
of skein modules of 3-manifolds}, in Knots 90 (Osaka, 1990), de Gruyter, 1992.

\item{8.} K. Johannson, {\it Homotopy Equivalences of
$3$--Manifolds with Boundaries}, L. N. M., vol. {\bf 761}, Springer
Verlag, 1979.

\item{9.} V.F.R. Jones, {\it A polynomial invariant for
knots via von Neumann algebras}, Bull. AMS, {\bf 12}(1985), 103--111.

\item{10.} $\underline {\phantom {aaaaaaaaaa}}$, {\it Hecke
algebra representations of braid groups and link polynomials},
Ann. of Math., {\bf 126}(1987), 335-338.

\item{11.} $\underline {\phantom {aaaaaaaaaa}}$, {\it On knot invariants
related to some statistical mechanical models},
Pacific J. of Math., {\bf 137}(1989), 311-334.

\item{12.} E. Kalfagianni, {\it Finite type invariants
for knots in 3-manifolds}, preprint, 1994.

\item{13.} M. Kontsevich, {\it Vassiliev's knot invariants},
Adv. in Sov. Math., {\bf 16}(1993), 137-150.

\item{14.} W.B.R. Lickorish, {\it A representation of orientable
combinatorial 3-manifolds}, Ann. of Math., {\bf 76}(1967), 531-538.

\item{15.} X.-S. Lin, {\it Finite
type link invariants of 3-manifolds}, Topology, {\bf 33}(1994), 45-71.

%\item{17.}  N.Y. Reshetikhin and
%V.G. Turaev,
%{ \it Invariants of $3$-manifolds via link polynomials and quantum groups},
%Invent. Math. {\bf 103} (1991), no. {\bf 3}, 547--597.

\item{16.} P. Scott, {\it There are no fake Seifert fibre spaces
with infinite $\pi_1$}, Ann. of Math., {\bf 117}(1983), 35-70.

\item{17.} V.G. Turaev, {\it The Yang-Baxter equation
and invariants of links}, Invent. Math., {\bf 92}(1988), 527--553.

\item{18.} V. A. Vassiliev, {\it Cohomology of Knot spaces},
in Theory of Singularities and its Applications (ed. V. Arnold), AMS, 1990.

%\item{22.} $\underline {\phantom {aaaaaaaaaa}}$, {\it
%Complements of Discriminants of Smooth Maps}, Topology and
%Application, AMS, 1992.

%\item {[W]} H. Wenzl, {\sl Representations of braid groups and the
%quantum Yang--Baxter equation}, Pac. J. Math. {\bf 145},
%No 1 (1990), 153-180.

\item{19.} E. Witten, {\it Quantum field theory
and Jones polynomial}, Commun. Math. Phys., {\bf 121} (1989),
359-389.
\bigskip

\p{\Smaller {\small S}CHOOL OF {\small M}ATHEMATICS,
 {\small I}NSTITUTE {\small F}OR
{\small A}DVANCED {\small S}TUDY {\small P}RINCETON, {\small NJ}
{\small 08540}}

\p{\small {\it Current address:}}

\p{\Smaller {\small D}EPARTMENT OF {\small M}ATHEMATICS,
 {\small R}UTGERS
{\small U}NIVERSITY, {\small N}EW {\small B}RUNSWICK, 
{\small NJ}
{\small 08903}}

\p{\small {E--mail address:}} {\vsmall  ekal@math.rutgers.edu}
\medskip

\p{\Smaller {\small D}EPARTMENT OF {\small M}ATHEMATICS,
 {\small C}OLUMBIA
{\small U}NIVERSITY {\small N}EW {\small Y}ORK, {\small NY}
{\small 10027}}

\p{\small {\it Current address:}}

\p{\Smaller {\small D}EPARTMENT OF {\small M}ATHEMATICS,
 {\small U}NIVERSITY OF
{\small C}ALIFORNIA, {\small R}IVERSIDE, 
{\small CA}
{\small 92521}}

\p{\small {E--mail address:}} {\vsmall  xl@math.columbia.edu}

\end

\end